# ACCURACY OF STATE SPACE COLLAPSE FOR EARLIEST-DEADLINE-FIRST QUEUES

By Lukasz Kruk,[1] John Lehoczky[2] and Steven Shreve[3]

*Maria Curie-Sklodowska University, Carnegie Mellon University and Carnegie Mellon University*

This paper presents a second-order heavy traffic analysis of a single server queue that processes customers having deadlines using the earliest-deadline-first scheduling policy. For such systems, referred to as *real-time queueing systems*, performance is measured by the fraction of customers who meet their deadline, rather than more traditional performance measures, such as customer delay, queue length or server utilization. To model such systems, one must keep track of customer lead times (the time remaining until a customer deadline elapses) or equivalent information. This paper reviews the earlier heavy traffic analysis of such systems that provided approximations to the system's behavior. The main result of this paper is the development of a second-order analysis that gives the accuracy of the approximations and the rate of convergence of the sequence of real-time queueing systems to its heavy traffic limit.

## 1. Introduction.

### 1.1. *Background.*

Real-time queueing systems are queueing systems whose customers have specific timing requirements. These systems arise in voice and video communication systems, control systems including avionics and automotive, and many aspects of modern manufacturing systems. The performance measures associated with such systems relate to the ability of the system to meet the customers' timing requirements as a function of the workload and the queue discipline. This is quite different from the more

Received October 2004; revised August 2005.

[1]Supported in part by the State Committee for Scientific Research of Poland, Grant 2 P03A 012 23.

[2]Supported in part by ONR and DARPA under MURI contract N00014-01-1-0576.

[3]Supported in part by NSF Grant DMS-04-04682.

common queueing system performance measures, such as customer delay, queue length or server utilization. In real-time queues, a queue discipline such as earliest deadline first (EDF) should be used rather than a more standard queue discipline like first-in-first-out (FIFO), which ignores the customer deadlines. EDF is a well studied queue discipline, especially in computer systems; see, for example, the monograph by Stankovic, Spuri, Ramamritham and Buttazzo [15]. In certain situations, EDF is optimal. For example, Panwar and Towsley [12] showed that, for a G/M/c queue with preemption, serving customers with general deadlines, EDF maximizes the fraction of customers that meet their deadlines within the class of work-conserving policies. In this paper we also treat the case of EDF with preemption.

In real-time queues, the state space is high dimensional in that the individual task deadlines, the task lead times (time until the deadline elapses) or some equivalent information must be kept to determine whether a customer deadline was met when it completes its service. This high-dimensional state space makes exact analysis intractable; however, a heavy traffic analysis can be carried out. This was done by Doytchinov, Lehoczky and Shreve (DLS) [4] for the single node, single traffic flow case, Kruk, Lehoczky, Shreve and Yeung [10] for the multi-class single station case, Yeung and Lehoczky [17] for feed-forward networks, and by Kruk, Lehoczky, Shreve and Yeung [11] for multi-flow acyclic networks. Both EDF and FIFO queue disciplines were considered.

In the rest of this section we summarize some of the existing results of the heavy traffic analysis and introduce the second-order analysis that is the subject of this paper. The model, assumptions and analysis are reviewed and made precise in Section 2 of this paper. The main results are stated in Section 3, and the proofs are developed in Sections 4–6. Simulations demonstrating the main results are presented in Section 7.

1.2. *Previous analytic results.* The heavy traffic analysis of real-time queues begins with a sequence of queueing systems, the $n$th system having independent strictly positive interarrival times with arrival rate $\lambda^{(n)}$, and the customers having independent service times with mean $\frac{1}{\mu^{(n)}}$. Assume the traffic intensity $\rho^{(n)} = \frac{\lambda^{(n)}}{\mu^{(n)}} = 1 - \frac{\gamma^{(n)}}{\sqrt{n}}$ for some sequence $\gamma^{(n)}$ having a limit $\gamma$. It follows that the scaled workload process

$$(1.1) \qquad \widehat{W}^{(n)}(t) \triangleq \frac{W(nt)}{\sqrt{n}} \Rightarrow W^*(t),$$

where $W^*$ is a reflected Brownian motion process with drift $-\gamma$ [see (2.17)].

To study the lead times of the customers in the queue, we introduce measure-valued processes $\mathcal{Q}^{(n)}(t)(B)$ and $\mathcal{W}^{(n)}(t)(B)$, where $\mathcal{Q}^{(n)}(t)(B)$ gives the number of customers in the queue at time $t$ having lead times in the Borel set $B$, while $\mathcal{W}^{(n)}(t)(B)$ is the work at time $t$ associated with customers in the queue having lead times in $B$. By considering $B = (-\infty, y], -\infty < y < \infty$, one can construct the cumulative lead-time distribution of work in the queue. The interval $(-\infty, 0)$ is of special importance because it is associated with work that is late.

To characterize the limiting behavior of these measure-valued processes, it is convenient to define the *frontier*, $F^{(n)}(t)$, roughly the largest lead time of all the customers ever having received any service. Any customer with lead time larger than $F^{(n)}(t)$ has never received any service. The frontier allows us to divide the workload (or customers) into two parts: those customers with lead times not larger than $F^{(n)}(t)$, that is, $\mathcal{W}^{(n)}(t)(-\infty, F^{(n)}(t)]$, and those with lead times greater than $F^{(n)}(t)$, that is, $\mathcal{W}^{(n)}(t)(F^{(n)}(t), \infty)$. DLS [4] prove that the scaled version of $\mathcal{W}^{(n)}(-\infty, F^{(n)}(t)]$, namely, $\frac{1}{\sqrt{n}}\mathcal{W}^{(n)}(nt)(-\infty, F^{(n)}(nt)]$, has limit zero; hence, we focus on the behavior of $\mathcal{W}^{(n)}(t)(y \vee F^{(n)}(t), \infty)$. (We use the notation $a \vee b \triangleq \max\{a, b\}$ and $a \wedge b \triangleq \min\{a, b\}$.)

Under the heavy traffic scaling implicit in the discussion so far, in order to obtain nontrivial limits as $n \to \infty$, time is accelerated by the factor $n$ and the workload and number of customers in queue is scaled by $\frac{1}{\sqrt{n}}$. Because the unscaled workload is of order $\sqrt{n}$, the time each customer spends in queue is also of order $\sqrt{n}$. In order to have a nontrivial limiting lead-time distribution, the lead times of arriving customers in the $n$th system must be of order $\sqrt{n}$. We assume therefore that arriving customers in the $n$th sytem have lead times equal to $\sqrt{n}$ times random variables drawn independently from a cumulative distribution function $G$. We assume that $G$ satisfies $y^* \triangleq \min\{y|G(y) = 1\} < \infty$, so that these random variables are bounded from above. It is natural to assume that $G(0-) = 0$, so that all lead times are nonnegative, but we do not need this for our analysis and, hence, do not assume it.

We set

$$(1.2) \quad H(y) \triangleq \int_y^\infty (1 - G(\eta))\, d\eta = \begin{cases} \displaystyle\int_y^{y^*} (1 - G(\eta))\, d\eta, & \text{if } y \leq y^*, \\ 0, & \text{if } y > y^*. \end{cases}$$

The function $H$ maps $(-\infty, y^*]$ onto $[0, \infty)$ and is strictly decreasing and Lipschitz continuous with Lipschitz constant 1 on $(-\infty, y^*]$. Therefore, there exists a continuous inverse function $H^{-1}$ that maps $[0, \infty)$ onto $(-\infty, y^*]$. In [4] it is shown that, as $n \to \infty$, the *scaled frontier process*

$$(1.3) \qquad \widehat{F}^{(n)}(t) \triangleq \frac{1}{\sqrt{n}} F^{(n)}(nt),$$

the *scaled workload measure process*

$$(1.4) \qquad \widehat{\mathcal{W}}^{(n)}(t)(B) \triangleq \frac{1}{\sqrt{n}} \mathcal{W}^{(n)}(nt)(\sqrt{n}B),$$

and the *scaled queue length measure process*

$$(1.5) \qquad \widehat{\mathcal{Q}}^{(n)}(t)(B) \triangleq \frac{1}{\sqrt{n}} \mathcal{Q}^{(n)}(nt)(\sqrt{n}B)$$

converge weakly to the *limiting scaled frontier process*

$$(1.6) \qquad F^*(t) \triangleq H^{-1}(W^*(t)), \qquad t \geq 0,$$

the *limiting scaled workload measure process*

$$(1.7) \qquad \widehat{\mathcal{W}}^*(t)(B) \triangleq \int_{B \cap [F^*(t), \infty)} (1 - G(y)) \, dy,$$

and the *limiting scaled queue-length measure process*

$$(1.8) \qquad \widehat{\mathcal{Q}}^*(t)(B) \triangleq \mu \int_{B \cap [F^*(t), \infty)} (1 - G(y)) \, dy,$$

where $\lambda = \lim_{n \to \infty} \lambda^{(n)}$ and $\mu = \lim_{n \to \infty} \mu^{(n)}$. Since $\lim_{n \to \infty} \rho^{(n)} = 1$, we have $\lambda = \mu$, but shall use both of these symbols to keep track of whether a term appears as the limit of $\lambda^{(n)}$ or the limit of $\mu^{(n)}$. This is useful to help conjecture the correct formulas if there were multiple input streams, in which case the limit of the sum of the arrival rates for the streams would equal $\mu$, rather than the limit for any particular arrival stream. It also aids in simulation (see Section 7), in which we are not yet at the limit and, thus, must replace $\lambda$ and $\mu$ in certain formulas by $\lambda^{(n)}$ and $\mu^{(n)}$, and these two are generally different.

In summary, DLS [4] prove that, as $n \to \infty$,

$$(1.9) \qquad \widehat{F}^{(n)} \Rightarrow F^*, \qquad \widehat{\mathcal{W}}^{(n)} \Rightarrow \widehat{\mathcal{W}}^* \quad \text{and} \quad \widehat{\mathcal{Q}}^{(n)} \Rightarrow \widehat{\mathcal{Q}}^*.$$

These convergence results allow us to approximate $\widehat{\mathcal{W}}^{(n)}(y, \infty)$ for $-\infty < y < \infty$ and $\widehat{F}^{(n)}$ in terms of the scaled workload process $\widehat{W}^{(n)}$. We present those approximations and discuss their accuracy below.

1.3. *Second-order analysis.* Our goal is to approximate the scaled frontier $\widehat{F}^{(n)}$ and the scaled workload measure $\widehat{\mathcal{W}}^{(n)}$ using the scaled workload scalar $\widehat{W}^{(n)}$, and to determine the accuracy of these approximations. The workload process is most useful as an approximation quantity as it is the most easily measured. Recall $\widehat{W}^{(n)}(t) \Rightarrow W^*(t)$, $\widehat{F}^{(n)}(t) \Rightarrow H^{-1}(W^*(t))$, and

$$\widehat{\mathcal{W}}^{(n)}(t)(y, \infty) \Rightarrow \widehat{\mathcal{W}}^*(t)(y, \infty) = H(y \vee F^*(t)) = H(y \vee H^{-1}(W^*(t))).$$

These suggest the following approximations for $\widehat{\mathcal{W}}^{(n)}(t)(y, \infty)$ and $\widehat{F}^{(n)}$:

$$(1.10) \qquad \widehat{\mathcal{W}}^{(n)}(t)(y, \infty) \approx H(y \vee \widehat{F}^{(n)}(t)),$$

$$(1.11) \qquad \widehat{\mathcal{W}}^{(n)}(t)(y, \infty) \approx H(y \vee H^{-1}(\widehat{W}^{(n)}(t))),$$

$$(1.12) \qquad \widehat{F}^{(n)}(t) \approx H^{-1}(\widehat{W}^{(n)}(t)).$$

In this paper we study the accuracy of these approximations by showing that the difference between the desired quantity and its proposed approximation when scaled by $n^{1/4}$ converges to a limiting process. As $n \to \infty$, the approximations in (1.10)–(1.12) become exact, and so some dilation of the difference of the two sides is necessary in order to obtain a nontrivial limit. In these approximations, $t$ is held fixed and the process parameter $y$ has been scaled by $\sqrt{n}$ [the $\sqrt{n}B$ term in (1.4)–(1.5)]. This is actually a scaling of lead times, and when time is scaled by $\sqrt{n}$, space must be scaled by $\frac{1}{n^{1/4}}$ in order to obtain a central limit result. However, space has already been scaled by $\frac{1}{\sqrt{n}}$ [the $\frac{1}{\sqrt{n}}$ term in (1.4) and (1.5)]. We must therefore multiply by $n^{1/4}$ to partially cancel this and obtain a scaling of $\frac{1}{n^{1/4}}$. Specifically, using [1] and [14], we prove the following results:

$$(1.13) \qquad n^{1/4}[\widehat{\mathcal{W}}^{(n)}(t)(y,\infty) - H(y \vee \widehat{F}^{(n)}(t))] \Rightarrow J^*(y \vee F^*(t)),$$

$$(1.14) \qquad n^{1/4}[\widehat{\mathcal{W}}^{(n)}(t)(y,\infty) - H(y \vee H^{-1}(\widehat{W}^{(n)}(t)))] \Rightarrow J^*(y)\mathbb{I}_{\{F^*(t) \leq y\}},$$

$$(1.15) \qquad n^{1/4}[H^{-1}(\widehat{W}^{(n)}(t)) - \widehat{F}^{(n)}(t)] \Rightarrow \frac{J^*(F^*(t))}{1 - G(F^*(t))},$$

where $J^*$ is a mean-zero Gaussian process with continuous paths and a covariance function defined in Theorem 3.3 below. In (1.13) and (1.14), both the left- and right-hand sides are processes in the parameter $y \leq y^*$ with $t > 0$ fixed; the convergence is weak convergence in $D(-\infty, y^*]$. In (1.15), $t > 0$ is again fixed and we have weak convergence of random variables.

## 2. The model, assumptions and notation.

2.1. *The basic model.* We now specify the model and its assumptions precisely. We have a sequence of single-station queueing systems, each serving one class of customers. The queueing systems are indexed by superscript $(n)$. The interarrival times for the customer arrival process are $\{u_j^{(n)}\}_{j=1}^{\infty}$, a sequence of strictly positive, independent, identically distributed random variables with common mean $\frac{1}{\lambda^{(n)}}$ and standard deviation $\alpha^{(n)}$. The service times are $\{v_j^{(n)}\}_{j=1}^{\infty}$, another sequence of positive, independent, identically distributed random variables with common mean $\frac{1}{\mu^{(n)}}$ and standard deviation $\beta^{(n)}$. We assume that each queue is empty at time zero.

We define the *customer arrival times*

$$(2.1) \qquad S_0^{(n)} \triangleq 0, \qquad S_k^{(n)} \triangleq \sum_{i=1}^{k} u_i^{(n)}, \qquad k \geq 1,$$

the *customer arrival process*

$$(2.2) \qquad A^{(n)}(t) \triangleq \max\{k; S_k^{(n)} \leq t\}, \qquad t \geq 0,$$

and the *work arrival process*

$$(2.3) \qquad V^{(n)}(t) \triangleq \sum_{j=1}^{\lfloor t \rfloor} v_j^{(n)}, \qquad t \geq 0.$$

The work that has arrived to the queue by time $t$ is then $V^{(n)}(A^{(n)}(t))$.

Each customer arrives with an initial lead time $L_j^{(n)}$, the time between the arrival time and the deadline for completion of service for that customer. These initial lead times are independent and identically distributed with distribution given by

$$(2.4) \qquad \mathbb{P}\{L_j^{(n)} \leq \sqrt{n}y\} = G(y),$$

where $G$ is a right-continuous cumulative distribution function. We define

$$(2.5) \qquad y^* \triangleq \min\{y \in \mathbb{R}; G(y) = 1\},$$

and assume that $y^*$ is finite. We assume that, for every $n$, the sequences $\{u_j^{(n)}\}_{j=1}^\infty$, $\{v_j^{(n)}\}_{j=1}^\infty$ and $\{L_j^{(n)}\}_{j=1}^\infty$ are mutually independent.

We assume that customers are served using the earliest-deadline-first (EDF) queue discipline, that is, the server always serves the customer with the shortest lead time. Preemption occurs when a customer more urgent than the customer in service arrives (we assume preempt-resume). There is no set up, switch-over or other type of overhead. Late customers (customers with negative lead times) stay in queue until served to completion.

The *netput process*

$$(2.6) \qquad N^{(n)}(t) \triangleq V^{(n)}(A^{(n)}(t)) - t$$

measures the amount of work in queue at time $t$, provided that the server is never idle up to time $t$. The *cumulative idleness process*

$$(2.7) \qquad I^{(n)}(t) \triangleq - \inf_{0 \leq s \leq t} N^{(n)}(s)$$

gives the amount of time the server is idle, and adding this to the netput process, we obtain the *workload process*

$$(2.8) \qquad W^{(n)}(t) = N^{(n)}(t) + I^{(n)}(t),$$

which records the amount of work in the queue, taking server idleness into account. All the above processes are independent of the queue service discipline, provided that the server is never idle when there are customers in the queue. However, the *queue length process* $Q^{(n)}(t)$, which is the number of customers in the queue at time $t$, depends on the queue discipline. All these processes are *right-continuous with left-hand limits* (*RCLL*).

2.2. *Heavy traffic assumptions.* We assume that the following limits exist and are all positive:

$$(2.9) \qquad \begin{aligned} \lim_{n\to\infty} \lambda^{(n)} = \lambda, & \qquad \lim_{n\to\infty} \mu^{(n)} = \mu, \\ \lim_{n\to\infty} \alpha^{(n)} = \alpha, & \qquad \lim_{n\to\infty} \beta^{(n)} = \beta. \end{aligned}$$

Define the *traffic intensity* $\rho^{(n)} \triangleq \frac{\lambda^{(n)}}{\mu^{(n)}}$. We make the *heavy traffic assumption*

$$(2.10) \qquad \lim_{n\to\infty} \sqrt{n}(1 - \rho^{(n)}) = \gamma$$

for some $\gamma \in \mathbb{R}$. This implies that $\lambda = \mu$, although, for reasons explained in Section 1.3, we shall use both symbols. We use the notation $\rho = \lambda/\mu$ in certain formulas, even though $\rho = 1$. We also impose the *modified Lindeberg condition* on the interarrival and service times: for every $c > 0$,

$$(2.11) \qquad \begin{aligned} \lim_{n\to\infty} \mathbb{E}[(u_j^{(n)} - (\lambda^{(n)})^{-1})^2 \mathbb{I}_{\{|u_j^{(n)} - (\lambda^{(n)})^{-1}| > cn^{1/4}\}}] \\ = \lim_{n\to\infty} \mathbb{E}[(v_j^{(n)} - (\mu^{(n)})^{-1})^2 \mathbb{I}_{\{|v_j^{(n)} - (\mu^{(n)})^{-1}| > cn^{1/4}\}}] = 0. \end{aligned}$$

This condition is satisfied, for example, if $\sup_{n\in\mathbb{N}} \mathbb{E}(u_j^{(n)})^{2+\delta} < \infty$ and $\sup_{n\in\mathbb{N}} \mathbb{E}(v_j^{(n)})^{2+\delta} < \infty$ for some $\delta > 0$. Clearly, (2.11) implies the usual Lindeberg condition on the interarrival and service times:

$$(2.12) \qquad \begin{aligned} \lim_{n\to\infty} \mathbb{E}[(u_j^{(n)} - (\lambda^{(n)})^{-1})^2 \mathbb{I}_{\{|u_j^{(n)} - (\lambda^{(n)})^{-1}| > c\sqrt{n}\}}] \\ = \lim_{n\to\infty} \mathbb{E}[(v_j^{(n)} - (\mu^{(n)})^{-1})^2 \mathbb{I}_{\{|v_j^{(n)} - (\mu^{(n)})^{-1}| > c\sqrt{n}\}}] = 0 \end{aligned}$$

for every $c > 0$. It can be shown that (2.12) does not, in general, imply (2.11).

We introduce the *heavy traffic scaling* for the idleness, workload and queue length processes

$$(2.13) \qquad \begin{aligned} \widehat{I}^{(n)}(t) &= \frac{1}{\sqrt{n}} I^{(n)}(nt), \\ \widehat{W}^{(n)}(t) &= \frac{1}{\sqrt{n}} W^{(n)}(nt), \\ \widehat{Q}^{(n)}(t) &= \frac{1}{\sqrt{n}} Q^{(n)}(nt), \end{aligned}$$

and the *centered heavy traffic scaling* for the arrival processes

$$\widehat{S}^{(n)}(t) = \frac{1}{\sqrt{n}} \sum_{j=1}^{\lfloor nt \rfloor} \left( u_j^{(n)} - \frac{1}{\lambda^{(n)}} \right),$$

$$\widehat{V}^{(n)}(t) = \frac{1}{\sqrt{n}} \sum_{j=1}^{\lfloor nt \rfloor} \left( v_j^{(n)} - \frac{1}{\mu^{(n)}} \right),$$

$$\widehat{A}^{(n)}(t) = \frac{1}{\sqrt{n}} [A^{(n)}(nt) - \lambda^{(n)} nt].$$

We define also

$$(2.14) \qquad \widehat{N}^{(n)}(t) = \frac{1}{\sqrt{n}} [V^{(n)}(A^{(n)}(nt)) - nt].$$

Note that $\widehat{W}^{(n)}(t) = \widehat{N}^{(n)}(t) + \widehat{I}^{(n)}(t)$.

Theorem 3.1 in [13] and Theorem 7.3.2 in [16] imply that

$$(2.15) \qquad (\widehat{S}^{(n)}, \widehat{A}^{(n)}) \Rightarrow (S^*, A^*),$$

where $A^*$ is a Brownian motion with zero drift and variance $\alpha^2 \lambda^3$ per unit time and

$$(2.16) \qquad S^*(\lambda t) = -\frac{1}{\lambda} A^*(t), \qquad t \geq 0.$$

It is a standard result [6] that

$$(2.17) \qquad (\widehat{N}^{(n)}, \widehat{I}^{(n)}, \widehat{W}^{(n)}) \Rightarrow (N^*, I^*, W^*),$$

where $N^*$ is a Brownian motion with variance $(\alpha^2 \rho^2 + \beta^2)\lambda$ per unit time and drift $-\gamma$, $I^*(t) \triangleq -\min_{0 \leq s \leq t} N^*(s)$, and $W^*(t) = N^*(t) + I^*(t)$. In other words, $W^*$ is a reflected Brownian motion with drift $-\gamma$, and $I^*$ causes the reflection.

Here and elsewhere, the symbol $\Rightarrow$ denotes weak convergence of measures on the space $D(T, S)$ of functions from a set $T$ (which is either a closed interval in $\mathbb{R}$ or a closed rectangle in $\mathbb{R}^2$, both possibly unbounded) to a Polish space $S$ that are right-continuous with left limits. If $S = \mathbb{R}$, we shall write simply $D(T)$. Throughout this paper, we shall use two topologies on $D(T, S)$. In almost all places, the Skorohod $J_1$ topology will be employed. This topology is convenient for dealing with weak convergence to continuous processes. The definition of the $J_1$ topology can be found, for example, in [3, 5, 16] for $T \subseteq \mathbb{R}$ and in [1] for the case of a rectangle $T$ in $\mathbb{R}^2$. In the sequel, whenever we consider weak convergence in $D(T, S)$, we always assume that this space is endowed with the $J_1$ topology unless explicitly stated otherwise. In particular, (2.15) and (2.17) hold in the $J_1$ topology. However, in Theorem 3.4, the (weaker) $M_1$ topology on a half-line $T$ will be used. We need to use the latter topology to establish stochastic-process limits with unmatched jumps in the limit process, for example, in the case of functions with asymptotically vanishing maximal jumps converging to an indicator function. See [16] for a definition of the $M_1$ topology and more details. We usually take $T = [0, \infty)$ and $S = \mathbb{R}^d$, with appropriate dimension $d$ [e.g., $d = 2$ in (2.15) and $d = 3$ in (2.17)], unless explicitly stated otherwise. For stochastic processes $A$ and $B$ with sample paths in $D(T, S)$, we shall write $A \overset{d}{=} B$ if $A$ and $B$ have the same distribution on $D(T, S)$.

2.3. *Measure-valued processes and frontiers.* To study whether tasks or customers meet their timing requirements, one must keep track of customer lead times, where the lead time is the time remaining until the deadline elapses, that is,

$$\text{lead time} \ = \ \text{deadline} \ - \ \text{current time}.$$

In this section we define a collection of measure-valued processes that will be useful in the analysis of the instantaneous lead-time profile of the customers. *Queue length measure*:

$$\mathcal{Q}^{(n)}(t)(B) \triangleq \left\{ \begin{matrix} \text{Number of customers in the queue at time} \\ t \text{ having lead times at time } t \text{ in } B \subseteq \mathbb{R} \end{matrix} \right\}.$$

*Workload measure*:

$$\mathcal{W}^{(n)}(t)(B) \triangleq \left\{ \begin{matrix} \text{Work in the queue at time } t \text{ associated with} \\ \text{customers having lead times at time } t \text{ in } B \subseteq \mathbb{R} \end{matrix} \right\}.$$

*Customer arrival measure*:

$$\mathcal{A}^{(n)}(t)(B) \triangleq \left\{ \begin{matrix} \text{Number of arrivals by time } t, \text{ whether} \\ \text{or not still in the system at time } t, \\ \text{having lead times at time } t \text{ in } B \subseteq \mathbb{R} \end{matrix} \right\}.$$

*Workload arrival measure*:

$$\mathcal{V}^{(n)}(t)(B) \triangleq \left\{ \begin{matrix} \text{Work associated with all arrivals by time } t, \\ \text{whether or not still in the system at time } t, \\ \text{having lead times at time } t \text{ in } B \subseteq \mathbb{R} \end{matrix} \right\}.$$

The following relationships easily follow:

(2.18) $\qquad Q^{(n)}(t) = \mathcal{Q}^{(n)}(t)(\mathbb{R}), \qquad\qquad W^{(n)}(t) = \mathcal{W}^{(n)}(t)(\mathbb{R}),$

(2.19) $\qquad A^{(n)}(t) = \mathcal{A}^{(n)}(t)(\mathbb{R}), \qquad V^{(n)}(A^{(n)}(t)) = \mathcal{V}^{(n)}(t)(\mathbb{R}),$

$$\mathcal{A}^{(n)}(t)(B) = \sum_{j=1}^{A^{(n)}(t)} \mathbb{I}_{\{L_j^{(n)} - (t - S_j^{(n)}) \in B\}}$$

(2.20)

$$= \sum_{j=1}^{\infty} \mathbb{I}_{\{S_j^{(n)} \in B + t - L_j^{(n)}, \, S_j^{(n)} \leq t\}},$$

$$\mathcal{V}^{(n)}(t)(B) = \sum_{j=1}^{A^{(n)}(t)} v_j^{(n)} \mathbb{I}_{\{L_j^{(n)} - (t - S_j^{(n)}) \in B\}}$$

(2.21)
$$= \sum_{j=1}^{\infty} v_j^{(n)} \mathbb{I}_{\{S_j^{(n)} \in B + t - L^{(n)}, S^{(n)} \leq t\}}.$$

To study the behavior of the EDF queue discipline, it is useful to keep track of the lead time of the customer currently in service and the largest lead time of all customers, whether present or departed, who have ever been in service. We define the *frontier*

$$F^{(n)}(t) \triangleq \left\{ \begin{array}{c} \text{Largest current lead time of all customers who have ever} \\ \text{been in service, whether still present or not, or } \sqrt{n}y^* - t, \\ \text{if this quantity is larger than the former one} \end{array} \right\},$$

and the *current lead time*

$$C^{(n)}(t) \triangleq \left\{ \begin{array}{c} \text{Lead time of the customer in service,} \\ \text{or } F^{(n)}(t) \text{ if the queue is empty} \end{array} \right\}.$$

Prior to arrival of the first customer, $F^{(n)}(t) = \sqrt{n}\,y^* - t$. Under the EDF queue discipline, there is no customer with lead time smaller than $C^{(n)}(t)$, and there has never been a customer in service whose lead time, if the customer were still present, would exceed $F^{(n)}(t)$. Furthermore, $C^{(n)}(t) \leq F^{(n)}(t)$ for all $t \geq 0$. Both $F^{(n)}$ and $C^{(n)}$ are RCLL.

We define the scaled versions (1.3), (1.4), (1.5) and $\widehat{C}^{(n)}(t) \triangleq \frac{1}{\sqrt{n}} C^{(n)}(nt)$ of the processes defined above under the EDF queue discpline. We define also

$$\widehat{\mathcal{A}}^{(n)}(t)(B) \triangleq \frac{1}{\sqrt{n}} \mathcal{A}^{(n)}(nt)(\sqrt{n}B)$$

(2.22)
$$= \frac{1}{\sqrt{n}} \sum_{j=1}^{A^{(n)}(nt)} \mathbb{I}_{\{L_j^{(n)} - (nt - S_j^{(n)}) \in \sqrt{n}B\}}$$

$$= \frac{1}{\sqrt{n}} \sum_{j=1}^{\infty} \mathbb{I}_{\{S_j^{(n)} \in \sqrt{n}B + nt - L_j^{(n)}, S_j^{(n)} \leq nt\}},$$

$$\widehat{\mathcal{V}}^{(n)}(t)(B) \triangleq \frac{1}{\sqrt{n}} \mathcal{V}^{(n)}(nt)(\sqrt{n}B)$$

(2.23)
$$= \frac{1}{\sqrt{n}} \sum_{j=1}^{A^{(n)}(nt)} v_j^{(n)} \mathbb{I}_{\{L_j^{(n)} - (nt - S_j^{(n)}) \in \sqrt{n}B\}}$$

$$= \frac{1}{\sqrt{n}} \sum_{j=1}^{\infty} v_j^{(n)} \mathbb{I}_{\{S_j^{(n)} \in \sqrt{n}B + nt - L_j^{(n)}, S_j^{(n)} \leq nt\}}.$$

**3. Main results.** Denote by $\mathcal{M}$ the set of all finite, nonnegative measures on $\mathcal{B}(\mathbb{R})$, the Borel subsets of $\mathbb{R}$. Under the weak topology, $\mathcal{M}$ is metrizable as a complete, separable topological space. We recall Proposition 3.10 and Theorem 3.1 of [4], which characterize the limiting distributions of the workload measure and the queue length measure under the EDF service discipline.

PROPOSITION 3.1 (Proposition 3.10 of [4]). *We have $\widehat{F}^{(n)} \Rightarrow F^*$ as $n \to \infty$.*

THEOREM 3.2 (Theorem 3.1 of [4]). *Let $\widehat{\mathcal{W}}^*$ and $\widehat{\mathcal{Q}}^*$ be the measure-valued processes defined by*

$$(3.1) \quad \widehat{\mathcal{W}}^*(t)(B) \triangleq \int_{B \cap [F^*(t), \infty)} (1 - G(y)) \, dy, \qquad \widehat{\mathcal{Q}}^*(t)(B) \triangleq \mu \widehat{\mathcal{W}}^*(t)(B),$$

*for all Borel sets $B \subseteq \mathbb{R}$. The processes $\widehat{\mathcal{W}}^{(n)}$ and $\widehat{\mathcal{Q}}^{(n)}$ converge weakly in $D([0, \infty), \mathcal{M})$ to $\widehat{\mathcal{W}}^*$ and $\widehat{\mathcal{Q}}^*$, respectively.*

By Theorem 3.2, for every $t \geq 0$ and $y \in \mathbb{R}$,

$$\widehat{\mathcal{W}}^{(n)}(t)(y, \infty) \Rightarrow \widehat{\mathcal{W}}^*(t)(y, \infty) = H(y \vee F^*(t)) = H(y \vee H^{-1}(W^*(t))).$$

In particular,

$$(3.2) \qquad\qquad \widehat{\mathcal{W}}^{(n)}(t)(y, \infty) \overset{d}{\approx} H(y \vee \widehat{F}^{(n)}(t)).$$

In some cases, the frontier $F^{(n)}$ (and, thus, its rescaled counterpart) may be difficult to evaluate. By Proposition 3.1, (1.6) and Theorem 3.2, we can replace $\widehat{F}^{(n)}(t)$ in (3.2) by the *approximate rescaled frontier* $H^{-1}(\widehat{W}^{(n)}(t))$, getting

$$(3.3) \qquad\qquad \widehat{\mathcal{W}}^{(n)}(t)(y, \infty) \approx H(y \vee H^{-1}(\widehat{W}^{(n)}(t))).$$

The aim of this paper is to investigate the rate of convergence of $\widehat{\mathcal{W}}^{(n)}(t)$ to $\widehat{\mathcal{W}}^*(t)$ in Theorem 3.2. More precisely, we find the empirical processes corresponding to the workload measure. In what follows, we fix $t > 0$. Our main results are the following:

THEOREM 3.3. *As a process in $y \leq y^*$, we have the convergence*

$$(3.4) \quad n^{1/4}[\widehat{\mathcal{W}}^{(n)}(t)(y, \infty) - H(y \vee \widehat{F}^{(n)}(t))] \Rightarrow J^*(y \vee F^*(t))$$

*in $D(-\infty, y^*]$, where $J^*$ is a mean-zero Gaussian process with continuous paths and covariance*

$$\mathbb{E}[J^*(y_1) J^*(y_2)] = \alpha^2 \lambda \rho^2 \int_{y_1}^{y^*} \int_{y_2}^{y^*} (\ell - y_1) \wedge (k - y_2) \, dG(\ell) \, dG(k)$$

$$(3.5) \qquad + \lambda \int_0^{y^* - (y_1 \vee y_2)} \left( \frac{1}{\mu^2} G((y_1 \wedge y_2) + x) + \beta^2 \right)$$
$$\times (1 - G((y_1 \vee y_2) + x)) \, dx,$$

*independent of $F^*(t)$.*

THEOREM 3.4. *As a process in $y \leq y^*$, we have the convergence*

$$(3.6) \quad n^{1/4}[\widehat{\mathcal{W}}^{(n)}(t)(y, \infty) - H(y \vee H^{-1}(\widehat{W}^{(n)}(t)))] \Rightarrow J^*(y)\mathbb{I}_{\{F^*(t) \leq y\}}$$

*in $D(-\infty, y^*]$ endowed with the $M_1$ topology, where $J^*$ is as in Theorem 3.3.*

Let us note that Theorems 3.3 and 3.4 characterize the accuracy of the approximations (3.2) and (3.3), respectively. One might also ask about the accuracy of approximating $\widehat{F}^{(n)}(t)$ by $H^{-1}(\widehat{W}^{(n)}(t))$. The answer to this question is given by the following proposition:

PROPOSITION 3.5. *For a fixed $t > 0$, we have the convergence*

$$(3.7) \qquad n^{1/4}[H^{-1}(\widehat{W}^{(n)}(t)) - \widehat{F}^{(n)}(t)] \Rightarrow \frac{J^*(F^*(t))}{1 - G(F^*(t))},$$

*where $J^*$ is as in Theorem 3.3.*

Let us observe that, by (1.6), $\mathbb{P}\{F^*(t) < y^*\} = \mathbb{P}\{W^*(t) > 0\} = 1$, so the limit in (3.7) is well defined.

**4. Customers behind the frontiers.** In this section we prove that the work in the $n$th system at time $nt$ associated with customers in this system having lead times smaller than the current frontier $F^{(n)}(nt)$ becomes negligible after division by $n^{1/4}$. The following lemma is a refinement of the second part of Proposition 3.6 in [4].

LEMMA 4.1. *Under the earliest-deadline-first queue discipline, we have*

$$(4.1) \qquad n^{1/4}\widehat{\mathcal{W}}^{(n)}(t)[\widehat{C}^{(n)}(t), \widehat{F}^{(n)}(t)) \Rightarrow 0.$$

PROOF. We define

$$(4.2) \qquad \tau^{(n)}(t) \triangleq \sup\{s \in [0, t]; \widehat{C}^{(n)}(s) = \widehat{F}^{(n)}(s)\}.$$

Because the system is initially empty, $\widehat{C}^{(n)}(0) = \widehat{F}^{(n)}(0) = 0$ and the above supremum is not over the empty set. In the proof of Proposition 3.6 in [4], it is shown that

$$(4.3) \qquad \sqrt{n}(t - \tau^{(n)}(t)) \Rightarrow 0.$$

Let us observe that, by the definition (4.2) and the fact that the interarrival times are strictly positive, we have

$$(4.4) \qquad \widehat{\mathcal{W}}^{(n)}(\tau^{(n)}(t))[\widehat{C}^{(n)}(\tau^{(n)}(t)), \widehat{F}^{(n)}(\tau^{(n)}(t))] \leq \frac{1}{\sqrt{n}} v^{(n)}_{A^{(n)}(n\tau^{(n)}(t))},$$

with strict inequality only when $\tau^{(n)}(t) = t$ and $\widehat{C}^{(n)}(t) = \widehat{F}^{(n)}(t)$, in which case the left-hand side is zero. We want to show that

$$(4.5) \qquad \frac{1}{n^{1/4}} v^{(n)}_{A^{(n)}(n\tau^{(n)}(t))} \Rightarrow 0.$$

To this end, let us choose a sequence $\varepsilon_n \downarrow 0$ so slowly that $P(A_n) \to 1$ as $n \to \infty$, where $A_n = [\tau^{(n)}(t) \geq t - \varepsilon_n/\sqrt{n}]$ [such a choice is possible by (4.3)]. By (2.15) and the differencing theorem (see, e.g., Appendix A of [4]),

$$\frac{1}{\sqrt{n}}[A^{(n)}(nt) - A^{(n)}(nt - \sqrt{n}\varepsilon_n) - \lambda^{(n)}\sqrt{n}\varepsilon_n] = \widehat{A}^{(n)}(t) - \widehat{A}^{(n)}\left(t - \frac{\varepsilon_n}{\sqrt{n}}\right) \Rightarrow 0.$$

Therefore, it is possible to find a sequence $a_n \downarrow 0$ so slowly that $P(B_n) \to 1$ as $n \to \infty$, where $B_n = [A^{(n)}(nt) - A^{(n)}(nt - \sqrt{n}\varepsilon_n) \leq b_n\sqrt{n}]$, and $b_n = \lambda^{(n)}\varepsilon_n + a_n$. For any $\delta > 0$, we have

$$
\begin{aligned}
&\lim_{n\to\infty} \mathbb{P}\left[\left\{\frac{1}{n^{1/4}} v^{(n)}_{A^{(n)}(n\tau^{(n)}(t))} \geq \delta\right\} \cap A_n \cap B_n\right] \\
&\qquad \leq \lim_{n\to\infty} \mathbb{P}\left[\left\{\frac{1}{n^{1/4}} \max_{A^{(n)}(nt) - b_n\sqrt{n} \leq j \leq A^{(n)}(nt)} v^{(n)}_j \geq \delta\right\} \cap A_n \cap B_n\right] \\
&\qquad \leq \lim_{n\to\infty} \mathbb{P}\left[\frac{1}{n^{1/4}} \max_{A^{(n)}(nt) - b_n\sqrt{n} \leq j \leq A^{(n)}(nt)} v^{(n)}_j \geq \delta\right] \\
(4.6) \quad &\qquad \leq \lim_{n\to\infty} \mathbb{P}\left[\frac{1}{n^{1/4}} \max_{1 \leq j \leq b_n\sqrt{n}+1} v^{(n)}_j \geq \delta\right] \\
&\qquad = 1 - \lim_{n\to\infty} \mathbb{P}\left[\bigcap_{j=1}^{\lfloor b_n\sqrt{n}\rfloor+1} \{v^{(n)}_j < \delta n^{1/4}\}\right] \\
&\qquad = 1 - \lim_{n\to\infty} \mathbb{P}[v^{(n)}_1 < \delta n^{1/4}]^{\lfloor b_n\sqrt{n}\rfloor+1} \\
&\qquad \leq 1 - \lim_{n\to\infty} \left(1 - \frac{\mathbb{E}(v^{(n)}_1)^2}{\delta^2\sqrt{n}}\right)^{\lfloor b_n\sqrt{n}\rfloor+1} = 0,
\end{aligned}
$$

where the fourth line follows from the independence of the service times and the arrivals, the sixth one from the fact that $\{v^{(n)}_j\}_{j=1}^{\infty}$ is a sequence of i.i.d. random variables, and the last one from (2.9) and the fact that $b_n \to 0$. But $P(A_n \cap B_n) \to 1$, so (4.6) implies (4.5).

An upper bound on the work with lead times in $[C^{(n)}(nt), F^{(n)}(nt))$ at time $nt$ is the work arrived to the system during the time interval $[n\tau^{(n)}(t), nt]$ minus the work served, which is $nt - n\tau^{(n)}(t)$. From this and (4.4), we obtain the bound

$$0 \leq n^{1/4}\widehat{\mathcal{W}}^{(n)}(t)[\widehat{C}^{(n)}(t), \widehat{F}^{(n)}(t))$$

$$\leq \frac{1}{n^{1/4}}v^{(n)}_{A^{(n)}(n\tau^{(n)}(t))} + \frac{1}{n^{1/4}}\sum_{j=A^{(n)}(n\tau^{(n)}(t))+1}^{A^{(n)}(nt)} v^{(n)}_j - n^{3/4}(t - \tau^{(n)}(t))$$

$$(4.7) \quad = \frac{1}{n^{1/4}}v^{(n)}_{A^{(n)}(n\tau^{(n)}(t))}$$

$$+ n^{1/4}\left[\widehat{V}^{(n)}\left(\frac{1}{n}A^{(n)}(nt)\right) - \widehat{V}^{(n)}\left(\frac{1}{n}A^{(n)}(n\tau^{(n)}(t))\right)\right]$$

$$+ \frac{1}{n^{1/4}\mu^{(n)}}[A^{(n)}(nt) - A^{(n)}(n\tau^{(n)}(t))] - n^{3/4}(t - \tau^{(n)}(t))$$

$$= \frac{1}{n^{1/4}}v^{(n)}_{A^{(n)}(n\tau^{(n)}(t))} + n^{1/4}\left[\widehat{V}^{(n)}\left(\frac{1}{\sqrt{n}}\widehat{A}^{(n)}(t) + \lambda^{(n)}t\right)\right.$$

$$\left. - \widehat{V}^{(n)}\left(\frac{1}{\sqrt{n}}\widehat{A}^{(n)}(\tau^{(n)}(t)) + \lambda^{(n)}\tau^{(n)}(t)\right)\right]$$

$$+ \frac{n^{1/4}}{\mu^{(n)}}[\widehat{A}^{(n)}(t) - \widehat{A}^{(n)}(\tau^{(n)}(t))] - (1 - \rho^{(n)})n^{3/4}(t - \tau^{(n)}(t)).$$

By (4.5), the first term on the right-hand side of (4.7) converges to zero in distribution. Also, by (2.10) and (4.3),

$$(4.8) \quad (1 - \rho^{(n)})n^{3/4}(t - \tau^{(n)}(t)) = n^{-1/4}(\sqrt{n}(1 - \rho^{(n)}))(\sqrt{n}(t - \tau^{(n)}(t)))$$

$$\Rightarrow 0.$$

Fix $\varepsilon > 0$. We have

$$n^{1/4}[\widehat{A}^{(n)}(t) - \widehat{A}^{(n)}(\tau^{(n)}(t))]$$

$$(4.9) \quad = n^{1/4}\left[\widehat{A}^{(n)}(t) - \widehat{A}^{(n)}\left(\left(t - \frac{\varepsilon}{\sqrt{n}}\right) \vee \tau^{(n)}(t)\right)\right]$$

$$+ n^{1/4}\left[\widehat{A}^{(n)}\left(\left(t - \frac{\varepsilon}{\sqrt{n}}\right) \vee \tau^{(n)}(t)\right) - \widehat{A}^{(n)}(\tau^{(n)}(t))\right].$$

The second term on the right-hand side of (4.9) converges weakly to zero by (4.3). The first term is bounded above by

$$n^{1/4} \max_{t - \varepsilon/\sqrt{n} \leq s \leq t} |\widehat{A}^{(n)}(t) - \widehat{A}^{(n)}(s)|$$

$$(4.10)$$

$$= \max_{0 \leq s \leq \varepsilon} n^{1/4}\left|\widehat{A}^{(n)}(t) - \widehat{A}^{(n)}\left(t - \frac{s}{\sqrt{n}}\right)\right|.$$

One can check that ordinary and renewal functional central limit theorems for triangular arrays (see, e.g., [5, 7, 13] and Theorem 14.6 in [3]) imply that

$$(4.11) \qquad \widetilde{A}^{(n)}(s) \triangleq n^{1/4}\left(\widehat{A}^{(n)}(t) - \widehat{A}^{(n)}\left(t - \frac{s}{\sqrt{n}}\right)\right) \Rightarrow B(s),$$

where $B$ is a Brownian motion (with zero drift and variance $\alpha^2\lambda^3$ per unit time). Therefore, (4.10) converges weakly to $\max_{0 \le s \le \varepsilon} |B(s)|$, which, in turn, converges to zero when $\varepsilon \downarrow 0$. We conclude that

$$(4.12) \qquad n^{1/4}[\widehat{A}^{(n)}(t) - \widehat{A}^{(n)}(\tau^{(n)}(t))] \Rightarrow 0.$$

The analysis of the second term on the right-hand side of (4.7) is similar to that given above. For a fixed $\varepsilon > 0$, we have

$$n^{1/4}\left[\widehat{V}^{(n)}\left(\frac{1}{\sqrt{n}}\widehat{A}^{(n)}(t) + \lambda^{(n)}t\right) - \widehat{V}^{(n)}\left(\frac{1}{\sqrt{n}}\widehat{A}^{(n)}(\tau^{(n)}(t)) + \lambda^{(n)}\tau^{(n)}(t)\right)\right]$$

$$= n^{1/4}\left[\widehat{V}^{(n)}\left(\frac{1}{\sqrt{n}}\widehat{A}^{(n)}(t) + \lambda^{(n)}t\right)\right.$$

$$- \widehat{V}^{(n)}\left(\frac{1}{\sqrt{n}}\widehat{A}^{(n)}\left(\left(t - \frac{\varepsilon}{\sqrt{n}}\right) \vee \tau^{(n)}(t)\right)\right.$$

$$(4.13) \qquad\qquad \left.\left. + \lambda^{(n)}\left(\left(t - \frac{\varepsilon}{\sqrt{n}}\right) \vee \tau^{(n)}(t)\right)\right)\right]$$

$$+ n^{1/4}\left[\widehat{V}^{(n)}\left(\frac{1}{\sqrt{n}}\widehat{A}^{(n)}\left(\left(t - \frac{\varepsilon}{\sqrt{n}}\right) \vee \tau^{(n)}(t)\right)\right.$$

$$\left. + \lambda^{(n)}\left(\left(t - \frac{\varepsilon}{\sqrt{n}}\right) \vee \tau^{(n)}(t)\right)\right)$$

$$\left. - \widehat{V}^{(n)}\left(\frac{1}{\sqrt{n}}\widehat{A}^{(n)}(\tau^{(n)}(t)) + \lambda^{(n)}\tau^{(n)}(t)\right)\right].$$

As before, the second of the two terms on the right-hand side of (4.13) converges weakly to zero by (4.3). To analyze the first one, we define a sequence of processes (depending on the parameters $t \ge 0$ and $\varepsilon > 0$) by

$$\widetilde{V}_{t,\varepsilon}^{(n)}(s) \triangleq n^{-1/4}\left[V^{(n)}\left(\sqrt{n}s + \lambda^{(n)}n\left(t - \frac{\varepsilon}{\sqrt{n}}\right)^+\right)\right.$$

$$\left. - V^{(n)}\left(\lambda^{(n)}n\left(t - \frac{\varepsilon}{\sqrt{n}}\right)^+\right) - \frac{\lfloor\sqrt{n}s\rfloor}{\mu^{(n)}}\right].$$

By the above-mentioned functional central limit theorems for triangular arrays,

$$(4.14) \qquad \widetilde{V}_{t,\varepsilon}^{(n)}(s) \Rightarrow \widetilde{B}(s),$$

where $\tilde{B}$ is a Brownian motion (with drift zero and variance $\beta^2$ per unit time). Moreover,

$$
\begin{aligned}
n^{1/4}\bigg[ &\widehat{V}^{(n)}\bigg( \frac{1}{\sqrt{n}}\widehat{A}^{(n)}(t) + \lambda^{(n)}t \bigg) \\
& - \widehat{V}^{(n)}\bigg( \frac{1}{\sqrt{n}}\widehat{A}^{(n)}\bigg( \bigg(t - \frac{\varepsilon}{\sqrt{n}}\bigg) \vee \tau^{(n)}(t) \bigg) \\
& \qquad + \lambda^{(n)}\bigg( \bigg(t - \frac{\varepsilon}{\sqrt{n}}\bigg) \vee \tau^{(n)}(t) \bigg) \bigg) \bigg]
\end{aligned}
$$

(4.15)
$$
\begin{aligned}
&= \widetilde{V}^{(n)}_{t,\varepsilon}\big( \widehat{A}^{(n)}(t) + \lambda^{(n)}(\sqrt{n}t - (\sqrt{n}t - \varepsilon)^+) \big) \\
& \quad - \widetilde{V}^{(n)}_{t,\varepsilon}\bigg( \widehat{A}^{(n)}\bigg( \bigg(t - \frac{\varepsilon}{\sqrt{n}}\bigg) \vee \tau^{(n)}(t) \bigg) \\
& \qquad\qquad + \lambda^{(n)}(\sqrt{n}\tau^{(n)}(t) - (\sqrt{n}t - \varepsilon)^+)^+ \bigg) \\
& \quad + O(n^{-1/4}).
\end{aligned}
$$

The right-hand side of (4.15) converges weakly to 0 by (2.15), (4.3), (4.14) and the differencing theorem. This, together with (4.7), (4.5), (4.8), (4.12) and (4.13), gives (4.1).  $\square$

## 5. Approximation for the workload arrival measure.
By Proposition 3.4 of [4], for every $y_0 < y^*$,

(5.1)
$$
\sup_{y_0 \le y \le y^*} |\widehat{\mathcal{V}}^{(n)}(t)(y,\infty) - H(y)| \xrightarrow{P} 0.
$$

In this section we want to find the joint limiting distribution for the the rescaled workload $\widehat{W}^{(n)}(t)$ and the *empirical process*

(5.2)
$$
\widehat{\mathcal{J}}^{(n)}(y) \triangleq n^{1/4}(\widehat{\mathcal{V}}^{(n)}(t)(y,\infty) - H(y)), \qquad y \le y^*,
$$

corresponding to (5.1). For $y \le y^*$, let

$$
M^{(n)}_j(y) \triangleq v^{(n)}_j \mathbb{I}_{\{L^{(n)}_j \le \sqrt{n}y\}} - \frac{1}{\mu^{(n)}}G(y).
$$

We have

(5.3)
$$
\begin{aligned}
\widehat{\mathcal{V}}^{(n)}(t)(y,\infty) &= \frac{1}{\sqrt{n}}\sum_{j=1}^{\infty} v^{(n)}_j \mathbb{I}_{\{nt+\sqrt{n}y-L^{(n)}_j < S^{(n)}_j \le nt\}} \\
&= I^{(n)}_1(y) + I^{(n)}_2(y),
\end{aligned}
$$

where

$$(5.4) \qquad I_1^{(n)}(y) = \int_y^{y^*} \frac{1}{\sqrt{n}} \sum_{j=1}^{\infty} \mathbb{I}_{\{nt+\sqrt{n}y-\sqrt{n}\ell < S_j^{(n)} \le nt\}} \, dM_j^{(n)}(\ell),$$

$$(5.5) \qquad I_2^{(n)}(y) = \int_y^{y^*} \frac{1}{\sqrt{n}} \sum_{j=1}^{\infty} \frac{1}{\mu^{(n)}} \mathbb{I}_{\{nt+\sqrt{n}y-\sqrt{n}\ell < S_j^{(n)} \le nt\}} \, dG(\ell).$$

In Sections 5.1 and 5.2 we study the limiting behavior of $I_1^{(n)}$ and $I_2^{(n)}$, respectively. We shall see that in the limit $n^{1/4} I_1^{(n)}$ depends only on the service times and lead times, not on the interarrival times, whereas $I_2^{(n)}$ obviously depends only on the interarrival times. Hence, the limits of $n^{1/4} I_1^{(n)}$ and $n^{1/4}(I_2^{(n)} - H)$, both of which exist, are independent. In Section 5.3 we use these obervations to characterize the limiting distribution of $(\widehat{J}^{(n)}, \widehat{W}^{(n)}(t))$.

5.1. *Asymptotic analysis for $I_1^{(n)}$.* For $\ell \ge y$, we have

$$A^{(n)}(nt) - A^{(n)}(nt - \sqrt{n}(\ell - y))$$

$$(5.6) \qquad = \lambda^{(n)} \sqrt{n}(\ell - y) + \sqrt{n}\left( \widehat{A}^{(n)}(t) - \widehat{A}^{(n)}\left( t - \frac{1}{\sqrt{n}}(\ell - y) \right) \right)$$

$$= \lambda^{(n)} \sqrt{n}(\ell - y) + n^{1/4} \widetilde{A}^{(n)}(\ell - y).$$

From (5.4) and (5.6), by the fact that customer arrival times are independent of their service times and lead times, we get

$$n^{1/4} I_1^{(n)}(y) = \frac{1}{n^{1/4}} \int_y^{y^*} \sum_{j=A^{(n)}(nt-\sqrt{n}(\ell-y))+1}^{A^{(n)}(nt)} dM_j^{(n)}(\ell)$$

$$(5.7) \qquad \stackrel{d}{=} \frac{1}{n^{1/4}} \int_y^{y^*} \sum_{j=1}^{A^{(n)}(nt)-A^{(n)}(nt-\sqrt{n}(\ell-y))} dM_j^{(n)}(\ell)$$

$$= \frac{1}{n^{1/4}} \int_y^{y^*} \sum_{j=1}^{\lambda^{(n)}\sqrt{n}(\ell-y)+n^{1/4}\widetilde{A}^{(n)}(\ell-y)} dM_j^{(n)}(\ell).$$

For $s \ge 0$ and $y \le y^*$, let us define a random field

$$(5.8) \quad Y^{(n)}(s,y) \triangleq \frac{1}{n^{1/4}} \sum_{j=1}^{\lfloor \lambda^{(n)} s \sqrt{n} \rfloor} \left( M_j^{(n)}(y^*) - M_j^{(n)}\left( \frac{j}{\lambda^{(n)}\sqrt{n}} + y \right) \right).$$

Then

$$\frac{1}{n^{1/4}} \int_y^{y^*} \sum_{j=1}^{\lfloor \lambda^{(n)}\sqrt{n}(\ell-y) \rfloor} dM_j^{(n)}(\ell)$$

$$(5.9) \quad \begin{aligned} &= \frac{1}{n^{1/4}} \int_y^{y^*} \sum_{j=1}^{\lfloor \lambda^{(n)}\sqrt{n}(y^*-y) \rfloor} \mathbb{I}_{\{\ell \geq j/(\lambda^{(n)}\sqrt{n})+y\}} \, dM_j^{(n)}(\ell) \\ &= \frac{1}{n^{1/4}} \sum_{j=1}^{\lfloor \lambda^{(n)}\sqrt{n}(y^*-y) \rfloor} \left( M_j^{(n)}(y^*) - M_j^{(n)}\left( \left( \frac{j}{\lambda^{(n)}\sqrt{n}} + y \right)- \right) \right) \\ &= Y^{(n)}(y^*-y, y-). \end{aligned}$$

Let

$$(5.10) \quad R^{(n)}(y) \triangleq \frac{1}{n^{1/4}} \int_y^{y^*} \left( \sum_{j=1}^{A^{(n)}(nt) - A^{(n)}(nt - \sqrt{n}(\ell-y))} - \sum_{j=1}^{\lfloor \lambda^{(n)}\sqrt{n}(\ell-y) \rfloor} \right) dM_j^{(n)}(\ell).$$

Then, by the second inequality in (5.7) and (5.8)–(5.10),

$$(5.11) \quad n^{1/4} I_1^{(n)}(y) \overset{d}{=} Y^{(n)}(y^*-y, y-) + R^{(n)}(y).$$

In the remainder of this subsection we find the limiting distribution for $Y^{(n)}$ and show that the process $R^{(n)}$ converges weakly to zero.

PROPOSITION 5.1. *For every $s_0 > 0$ and $y_0 < y^*$, $Y^{(n)} \Rightarrow Y$ in $D([0, s_0] \times [y_0, y^*])$, where $Y$ is a mean-zero Gaussian random field with continuous paths and covariance*

$$\begin{aligned} &\mathbb{E}[Y(s_1, y_1) Y(s_2, y_2)] \\ &\quad = \lambda \int_0^{s_1 \wedge s_2} \left( \frac{1}{\mu^2} G((y_1 \wedge y_2) + x) + \beta^2 \right) (1 - G((y_1 \vee y_2) + x)) \, dx. \end{aligned}$$

PROOF. We will first show that the sequence $\{Y^{(n)}\}$ is tight. For $y \leq y^*$ and $j = 1, 2, \dots$, let

$$G_j(y) \triangleq G\left( \frac{j}{\lambda^{(n)}\sqrt{n}} + y \right).$$

Also, for $0 \leq s_1 \leq s_2 \leq s_3 \leq s_0$, $y_0 \leq y_1 \leq y_2 \leq y_3 \leq y^*$, let

$$\begin{aligned} &Y^{(n)}((s_i, s_{i+1}) \times (y_k, y_{k+1})) \\ &\quad \triangleq Y^{(n)}(s_{i+1}, y_{k+1}) - Y^{(n)}(s_i, y_{k+1}) - Y^{(n)}(s_{i+1}, y_k) + Y^{(n)}(s_i, y_k) \\ &\quad = \frac{1}{n^{1/4}} \sum_{j=\lfloor \lambda^{(n)} s_i \sqrt{n} \rfloor + 1}^{\lfloor \lambda^{(n)} s_{i+1}\sqrt{n} \rfloor} \Delta M_j^{(n)}(y_k, y_{k+1}), \end{aligned}$$

where

$$\Delta M_j^{(n)}(y_k, y_{k+1}) \triangleq M_j^{(n)}\left( \frac{j}{\lambda^{(n)}\sqrt{n}} + y_k \right) - M_j^{(n)}\left( \frac{j}{\lambda^{(n)}\sqrt{n}} + y_{k+1} \right).$$

We define $B_1 = (s_1, s_2] \times (y_1, y_2]$ and we will take $B_2$ to be a block "neighboring" $B_1$, considering both the case that

$$(5.12) \qquad\qquad B_2 = (s_2, s_3] \times (y_1, y_2]$$

is to the right of $B_1$ and also the case that

$$(5.13) \qquad\qquad B_2 = (s_1, s_2] \times (y_2, y_3]$$

is above $B_1$. Our goal in both cases is to obtain the bound (5.23).

We continue with $B_2$ given by (5.12). The independence of the random variables $\Delta M_j^{(n)}(y_k, y_{k+1})$ for different values of $j$ implies

$$(5.14) \qquad \mathbb{E}(Y^{(n)}(B_i))^2 = \frac{1}{\sqrt{n}} \sum_{j=\lfloor \lambda^{(n)} s_i \sqrt{n} \rfloor + 1}^{\lfloor \lambda^{(n)} s_{i+1} \sqrt{n} \rfloor} \mathbb{E}(\Delta M_j^{(n)}(y_1, y_2))^2,$$

$$(5.15) \quad \mathbb{E}[(Y^{(n)}(B_1))^2 (Y^{(n)}(B_2))^2] = \mathbb{E}(Y^{(n)}(B_1))^2 \mathbb{E}(Y^{(n)}(B_2))^2,$$

for $i = 1, 2$. Let us note that

$$
\begin{aligned}
& \mathbb{E}(\Delta M_j^{(n)}(y_1, y_2))^2 \\
(5.16) \quad & = \mathbb{E}\left[ v_j^{(n)} \mathbb{1}_{\{j/\lambda^{(n)} + \sqrt{n} y_1 < L_j^{(n)} \leq j/\lambda^{(n)} + \sqrt{n} y_2\}} - \frac{1}{\mu^{(n)}} (G_j(y_2) - G_j(y_1)) \right]^2 \\
& = \mathbb{E}(v_j^{(n)})^2 (G_j(y_2) - G_j(y_1)) - \frac{1}{(\mu^{(n)})^2} (G_j(y_2) - G_j(y_1))^2 \\
& \leq C(G_j(y_2) - G_j(y_1)),
\end{aligned}
$$

where

$$C = \sup_n \mathbb{E}(v_j^{(n)})^2 = \sup_n \left\{ (\beta^{(n)})^2 + \frac{1}{(\mu^{(n)})^2} \right\} < \infty$$

is a constant independent of $j$ and $n$ because of (2.9). Let $U^{(n)}$ be a random variable such that

$$(5.17) \quad \mathbb{P}\left[ U^{(n)} = \frac{j}{\lambda^{(n)} \sqrt{n}} \right] = \frac{1}{\lfloor \lambda^{(n)} \sqrt{n} s_0 \rfloor}, \qquad j = 1, \ldots, \lfloor \lambda^{(n)} \sqrt{n} s_0 \rfloor.$$

Let $X$ be a random variable with cumulative distribution function $G$, independent of $U^{(n)}$. For $0 \leq s_1 \leq s_2 \leq s_0$, $y_0 \leq y_1 \leq y_2 \leq y^*$, let us define

$$
\begin{aligned}
& m_1^{(n)}((s_1, s_2] \times (y_1, y_2]) \\
(5.18) \quad & \triangleq \mathbb{P}[(U^{(n)}, X - U^{(n)}) \in (s_1, s_2] \times (y_1, y_2]] \\
& = \frac{1}{\lfloor \lambda^{(n)} \sqrt{n} s_0 \rfloor} \sum_{j=\lfloor s_1 \lambda^{(n)} \sqrt{n} \rfloor + 1}^{\lfloor s_2 \lambda^{(n)} \sqrt{n} \rfloor} (G_j(y_2) - G_j(y_1)).
\end{aligned}
$$

We have, by (5.14)–(5.18),

$$\mathbb{E}[(Y^{(n)}(B_1))^2(Y^{(n)}(B_2))^2] \leq C^2 \left(\frac{\lfloor \lambda^{(n)}\sqrt{n}s_0\rfloor}{\sqrt{n}}\right)^2 m_1^{(n)}(B_1)m_1^{(n)}(B_2)$$

(5.19)
$$\leq C^2 (\lambda^{(n)}s_0)^2 m_1^{(n)}(B_1)m_1^{(n)}(B_2).$$

We have obtained (5.23) with $m^{(n)} \triangleq \lambda^{(n)}s_0 C m_1^{(n)}$.

Now let $B_2$ be given by (5.13). Using the fact that $\mathbb{E}\Delta M_j^{(n)}(y_{k+1}, y_k) = 0$ and the $\Delta M_j^{(n)}(y_k, y_{k+1})$ are independent for different values of $j$, we have

$$\mathbb{E}[(Y^{(n)}(B_1))^2(Y^{(n)}(B_2))^2]$$

$$= \frac{1}{n}\mathbb{E}\left[\left(\sum_{j=\lfloor \lambda^{(n)}s_1\sqrt{n}\rfloor+1}^{\lfloor \lambda^{(n)}s_2\sqrt{n}\rfloor} \Delta M_j^{(n)}(y_1, y_2)\right)^2 \left(\sum_{j=\lfloor \lambda^{(n)}s_1\sqrt{n}\rfloor+1}^{\lfloor \lambda^{(n)}s_2\sqrt{n}\rfloor} \Delta M_j^{(n)}(y_2, y_3)\right)^2\right]$$

$$= \frac{1}{n}\sum_{j=\lfloor \lambda^{(n)}s_1\sqrt{n}\rfloor+1}^{\lfloor \lambda^{(n)}s_2\sqrt{n}\rfloor} \mathbb{E}[(\Delta M_j^{(n)}(y_1, y_2))^2(\Delta M_j^{(n)}(y_2, y_3))^2]$$

(5.20)
$$+ \frac{1}{n}\sum_{\substack{j,k=\lfloor \lambda^{(n)}s_1\sqrt{n}\rfloor+1 \\ j\neq k}}^{\lfloor \lambda^{(n)}s_2\sqrt{n}\rfloor} [\mathbb{E}(\Delta M_j^{(n)}(y_1, y_2))^2 \mathbb{E}(\Delta M_k^{(n)}(y_2, y_3))^2]$$

$$+ \frac{1}{n}\sum_{\substack{j,k=\lfloor \lambda^{(n)}s_1\sqrt{n}\rfloor+1 \\ j\neq k}}^{\lfloor \lambda^{(n)}s_2\sqrt{n}\rfloor} \mathbb{E}[\Delta M_j^{(n)}(y_1, y_2)\Delta M_j^{(n)}(y_2, y_3)]$$

$$\times \mathbb{E}[\Delta M_k^{(n)}(y_1, y_2)\Delta M_k^{(n)}(y_2, y_3)].$$

Using (5.20) and proceeding as in (5.16), we can check that there exists a constant $C_1$ independent of $j$ and $n$ such that

$$\mathbb{E}[(Y^{(n)}(B_1))^2(Y^{(n)}(B_2))^2]$$

(5.21)
$$\leq \frac{C_1}{n}\left\{\sum_{j=\lfloor \lambda^{(n)}s_1\sqrt{n}\rfloor+1}^{\lfloor \lambda^{(n)}s_2\sqrt{n}\rfloor} (G_j(y_2)-G_j(y_1))(G_j(y_3)-G_j(y_2))\right.$$

$$\left. + 2\sum_{\substack{j,k=\lfloor \lambda^{(n)}s_1\sqrt{n}\rfloor+1 \\ j\neq k}}^{\lfloor \lambda^{(n)}s_2\sqrt{n}\rfloor} (G_j(y_2)-G_j(y_1))(G_k(y_3)-G_k(y_2))\right\}.$$

Thus,

(5.22)   $\mathbb{E}[(Y^{(n)}(B_1))^2(Y^{(n)}(B_2))^2] \leq 2(s_0\lambda^{(n)})^2 C_1 m_1^{(n)}(B_1)m_1^{(n)}(B_2).$

By (5.18), (5.19) and (5.22), for all "neighboring blocks" $B_1$ and $B_2$ in $[0, s_0] \times [y_0, y^*]$, we have

$$(5.23) \qquad \mathbb{E}[(Y^{(n)}(B_1))^2 (Y^{(n)}(B_2))^2] \le m^{(n)}(B_1) m^{(n)}(B_2),$$

where $m^{(n)} \triangleq \lambda^{(n)} s_0 (C + \sqrt{2C_1}) m_1^{(n)}$.

It is clear that $m_1^{(n)} \Rightarrow m_1$, and thus,

$$(5.24) \qquad m^{(n)} \Rightarrow m,$$

where

$$(5.25) \quad \begin{aligned} m_1((s_1, s_2] \times (y_1, y_2]) &\triangleq \mathbb{P}[(U, X - U) \in (s_1, s_2] \times (y_1, y_2]], \\ m &\triangleq \lambda s_0 (C + \sqrt{2C_1}) m_1, \end{aligned}$$

with independent random variables $U$ and $X$ having distributions uniform on $[0, s_0]$ and $G$, respectively. In particular, $m$ has continuous marginals. Then, by (5.23), (5.24) and Theorem 3 from [1] (strictly speaking, by its extension described on pages 1665–1666 of that paper), the sequence $\{Y^{(n)}\}$ is tight in $D([0, s_0] \times [y_0, y^*])$. (Note that $Y^{(n)}$ vanishes along $(\{0\} \times [y_0, y^*]) \cup ([0, s_0] \times \{y^*\})$, instead of the "lower boundary" $(\{0\} \times [y_0, y^*]) \cup ([0, s_0] \times \{y_0\})$, as required by the assumptions of (the extension of) Theorem 3 from [1], but the proof of the latter result clearly goes through also in our case.)

Now we will show that the finite-dimensional distributions of $Y^{(n)}$ converge to the corresponding distributions of $Y$. Let $0 \le s_1 \le s_2 \le s_0$ and $y_0 \le y_1, y_2 \le y^*$ be given. We claim that

$$(5.26) \qquad \lim_{n \to \infty} \mathbb{E}[Y^{(n)}(s_1, y_1) Y^{(n)}(s_2, y_2)] = \mathbb{E}[Y(s_1, y_1) Y(s_2, y_2)].$$

Let us observe that

$$\begin{aligned} \mathbb{E}[Y^{(n)}(s_1, y_1) Y^{(n)}(s_2, y_2)] &= \mathbb{E}[Y^{(n)}(s_1, y_1) Y^{(n)}(s_1, y_2)] \\ &\quad + \mathbb{E}[Y^{(n)}(s_1, y_1)(Y^{(n)}(s_2, y_2) - Y^{(n)}(s_1, y_2))] \\ &= \mathbb{E}[Y^{(n)}(s_1, y_1) Y^{(n)}(s_1, y_2)], \end{aligned}$$

and $\mathbb{E}[Y(s_1, y_1) Y(s_2, y_2)] = \mathbb{E}[Y(s_1, y_1) Y(s_1, y_2)]$ by the definition of $Y$, so it suffices to check (5.26) for $s_1 = s_2 \triangleq s$. Assume that $y_1 \le y_2$. Then

$$\begin{aligned} &\mathbb{E}[Y^{(n)}(s, y_1) Y^{(n)}(s, y_2)] \\ &\quad = \frac{1}{\sqrt{n}} \sum_{j=1}^{\lfloor \lambda^{(n)} s \sqrt{n} \rfloor} \mathbb{E}\Bigg[ \bigg( M_j^{(n)}(y^*) - M_j^{(n)}\bigg( \frac{j}{\lambda^{(n)} \sqrt{n}} + y_1 \bigg) \bigg) \\ &\qquad\qquad \times \bigg( M_j^{(n)}(y^*) - M_j^{(n)}\bigg( \frac{j}{\lambda^{(n)} \sqrt{n}} + y_2 \bigg) \bigg) \Bigg] \end{aligned}$$

$$= \frac{1}{\sqrt{n}} \sum_{j=1}^{\lfloor \lambda^{(n)} s \sqrt{n} \rfloor} \mathbb{E}\Bigg[\Bigg( v_j^{(n)} \mathbb{I}_{\{j/(\lambda^{(n)}\sqrt{n})+y_1 < L_j^{(n)}/\sqrt{n}\}} - \frac{1}{\mu^{(n)}}(1 - G_j(y_1)) \Bigg)$$

$$\times \Bigg( v_j^{(n)} \mathbb{I}_{\{j/(\lambda^{(n)}\sqrt{n})+y_2 < L_j^{(n)}/\sqrt{n}\}} - \frac{1}{\mu^{(n)}}(1 - G_j(y_2)) \Bigg) \Bigg]$$

$$= \frac{1}{\sqrt{n}} \sum_{j=1}^{\lfloor \lambda^{(n)} s \sqrt{n} \rfloor} \Bigg[ \mathbb{E}(v_j^{(n)})^2 (1 - G_j(y_2))$$

$$- \frac{1}{(\mu^{(n)})^2}(1 - G_j(y_1))(1 - G_j(y_2)) \Bigg]$$

$$= \frac{1}{\sqrt{n}} \sum_{j=1}^{\lfloor \lambda^{(n)} s \sqrt{n} \rfloor} \Bigg[ \Bigg( \frac{1}{(\mu^{(n)})^2} + (\beta^{(n)})^2 \Bigg) - \frac{1}{(\mu^{(n)})^2}(1 - G_j(y_1)) \Bigg](1 - G_j(y_2))$$

$$= \frac{1}{\sqrt{n}} \sum_{j=1}^{\lfloor \lambda^{(n)} s \sqrt{n} \rfloor} \Bigg[ \frac{1}{(\mu^{(n)})^2} G_j(y_1) + (\beta^{(n)})^2 \Bigg](1 - G_j(y_2)).$$

To obtain (5.26), we observe that the last term in (5.27) converges to

$$\lambda \int_0^s \Bigg( \frac{1}{\mu^2} G(y_1 + x) + \beta^2 \Bigg)(1 - G(y_2 + x)) \, dx = \mathbb{E}[Y(s, y_1)Y(s, y_2)].$$

To show convergence of the finite-dimensional distributions of $Y^{(n)}$, we will use the Cramér–Wold device (see, e.g., [2]). Fix $m$, $0 \le s_1 \le \cdots \le s_m \le s_0$, $y_0 \le y_1, \ldots, y_m \le y^*$ and $t_1, \ldots, t_m \in \mathbb{R}$. Then

$$\sum_{i=1}^m t_i Y^{(n)}(s_i, y_i)$$

$$= \frac{1}{n^{1/4}} \sum_{i=1}^m t_i \sum_{j=1}^{\lfloor \lambda^{(n)} s_i \sqrt{n} \rfloor} \Bigg( M_j^{(n)}(y^*) - M_j^{(n)}\Bigg( \frac{j}{\lambda^{(n)}\sqrt{n}} + y_i \Bigg) \Bigg)$$

$$= \sum_{j=1}^{\lfloor \lambda^{(n)} s_m \sqrt{n} \rfloor} \tilde{X}_j^{(n)},$$

where

$$\tilde{X}_j^{(n)} \triangleq \frac{1}{n^{1/4}} \sum_{\{i\,:\,j \le \lambda^{(n)}\sqrt{n} s_i\}} t_i \Bigg( M_j^{(n)}(y^*) - M_j^{(n)}\Bigg( \frac{j}{\lambda^{(n)}\sqrt{n}} + y_i \Bigg) \Bigg)$$

are independent, mean-zero random variables. By (5.28), we have

$$\mathbf{s}_n^2 \triangleq \mathrm{Var}\Bigg( \sum_{j=1}^{\lfloor \lambda^{(n)} s_m \sqrt{n} \rfloor} \tilde{X}_j^{(n)} \Bigg) = \mathrm{Var}\Bigg( \sum_{i=1}^m t_i Y^{(n)}(s_i, y_i) \Bigg)$$

(5.29)
$$= \sum_{i,k=1}^{m} t_i t_k \mathbb{E}[Y^{(n)}(s_i, y_i) Y^{(n)}(s_k, y_k)],$$

so, by (5.26),

$$\mathbf{s}^2 \triangleq \lim_{n \to \infty} \mathbf{s}_n^2 = \sum_{i,k=1}^{m} t_i t_k \mathbb{E}[Y(s_i, y_i) Y(s_k, y_k)]$$

(5.30)
$$= \text{Var}\left( \sum_{i=1}^{m} t_i Y(s_i, y_i) \right).$$

If $\mathbf{s} = 0$, then $\sum_{i=1}^{m} t_i Y^{(n)}(s_i, y_i) \Rightarrow \sum_{i=1}^{m} t_i Y(s_i, y_i)$ by (5.29), (5.30) and the Chebyshev inequality. Assume $\mathbf{s} > 0$. We shall check that the random variables $\tilde{X}_j^{(n)}$ satisfy the Lindeberg condition. Let $X_j^{(n)} \triangleq n^{1/4} \tilde{X}_j^{(n)}$. It is easy to see that there exist constants $C_1, C_2 > 0$ such that $|X_j^{(n)}| \leq C_1 v_j^{(n)} + C_2$ for every $j$, $n$. Let $\varepsilon > 0$ be arbitrary. Then

$$\frac{1}{\mathbf{s}_n^2} \sum_{j=1}^{\lfloor \lambda^{(n)} s_m \sqrt{n} \rfloor} \int_{\{|\tilde{X}_j^{(n)}| \geq \varepsilon \mathbf{s}_n\}} (\tilde{X}_j^{(n)})^2 \, d\mathbb{P}$$

$$= \frac{1}{\mathbf{s}_n^2 \sqrt{n}} \sum_{j=1}^{\lfloor \lambda^{(n)} s_m \sqrt{n} \rfloor} \int_{\{|X_j^{(n)}| \geq \varepsilon \mathbf{s}_n n^{1/4}\}} (X_j^{(n)})^2 \, d\mathbb{P}$$

(5.31)
$$\leq \frac{2}{\mathbf{s}_n^2 \sqrt{n}} \sum_{j=1}^{\lfloor \lambda^{(n)} s_m \sqrt{n} \rfloor} \int_{\{v_j^{(n)} \geq (1/C_1)(\varepsilon \mathbf{s}_n n^{1/4} - C_2)\}} (C_1^2 (v_j^{(n)})^2 + C_2^2) \, d\mathbb{P}$$

$$= O(1) \int_{\{v_1^{(n)} \geq (1/C_1)(\varepsilon \mathbf{s}_n n^{1/4} - C_2)\}} (v_1^{(n)})^2 \, d\mathbb{P}$$

$$+ O(1) \mathbb{P}\left[ v_1^{(n)} \geq \frac{\varepsilon \mathbf{s}_n n^{1/4} - C_2}{C_1} \right].$$

As $n \to \infty$,

(5.32)
$$\mathbb{P}\left[ v_1^{(n)} \geq \frac{\varepsilon \mathbf{s}_n n^{1/4} - C_2}{C_1} \right] \leq \frac{C_1}{\varepsilon \mathbf{s}_n n^{1/4} - C_2} \mathbb{E} v_1^{(n)}$$

$$= \frac{C_1}{(\varepsilon \mathbf{s}_n n^{1/4} - C_2) \mu^{(n)}} \to 0$$

and

(5.33)
$$\int_{\{v_1^{(n)} \geq (1/C_1)(\varepsilon \mathbf{s}_n n^{1/4} - C_2)\}} (v_1^{(n)})^2 \, d\mathbb{P} \to 0$$

by (2.11) and (5.32). Thus, by (5.31)–(5.33), the random variables $\tilde{X}_j^{(n)}$ satisfy the Lindeberg condition. Therefore, by the Lindeberg central limit theorem, the sum $\sum_{i=1}^{m} t_i Y^{(n)}(s_i, y_i)$ converges weakly to $N(0, \mathbf{s})$, the distribution of $\sum_{i=1}^{m} t_i Y(s_i, y_i)$. We have proved convergence of the finite-dimensional distributions of $Y^{(n)}$ to the corresponding distributions of $Y$.

Finally, we show that the limiting random field $Y$ has continuous sample paths. By the Kolmogorov–Čentsov theorem (see, e.g., [8]), it suffices to show that there exists a constant $C$ such that, for any $0 \leq s_1 \leq s_2$, $y_1, y_2 \leq y^*$,

$$(5.34) \qquad \mathbb{E}|Y(s_1, y_1) - Y(s_2, y_2)|^6 \leq C \|(s_1, y_1) - (s_2, y_2)\|^3,$$

where $\|\cdot\|$ denotes the Euclidean norm in $\mathbb{R}^2$. It is well known that, for every $n$, there exists a constant $C_n$ such that, for every normal random variable $\tilde{Z}$ with mean zero, $\mathbb{E}\tilde{Z}^{2n} \leq C_n(\mathbb{E}\tilde{Z}^2)^n$. In particular, because $Y$ is Gaussian, to prove (5.34), it suffices to find a constant $\tilde{C}$ such that

$$(5.35) \qquad \mathbb{E}|Y(s_1, y_1) - Y(s_2, y_2)|^2 \leq \tilde{C} \|(s_1, y_1) - (s_2, y_2)\|.$$

We have

$$
\begin{aligned}
(5.36) \quad & \mathbb{E}|Y(s_1, y_1) - Y(s_2, y_2)|^2 \\
& = \mathbb{E}Y(s_1, y_1)^2 - 2\mathbb{E}[Y(s_1, y_1)Y(s_2, y_2)] + \mathbb{E}Y(s_2, y_2)^2.
\end{aligned}
$$

But

$$
\begin{aligned}
(5.37) \quad & |\mathbb{E}Y(s_2, y_2)^2 - \mathbb{E}[Y(s_1, y_1)Y(s_2, y_2)]| \\
& = \lambda \Bigg| \int_0^{s_2} \left( \frac{1}{\mu^2} G(y_2 + x) + \beta^2 \right)(1 - G(y_2 + x))\, dx \\
& \qquad\qquad - \int_0^{s_1} \left( \frac{1}{\mu^2} G((y_1 \wedge y_2) + x) + \beta^2 \right) \\
& \qquad\qquad\qquad\qquad \times (1 - G((y_1 \vee y_2) + x))\, dx \Bigg|.
\end{aligned}
$$

If $y_1 \leq y_2$, then the right-hand side of (5.37) equals

$$
\begin{aligned}
& \frac{\lambda}{\mu^2} \int_0^{s_1} (G(y_2 + x) - G(y_1 + x))(1 - G(y_2 + x))\, dx \\
& \quad + \lambda \int_{s_1}^{s_2} \left( \frac{1}{\mu^2} G(y_2 + x) + \beta^2 \right)(1 - G(y_2 + x))\, dx \\
& \quad \leq \frac{\lambda}{\mu^2} \int_0^{s_1} (G(y_2 + x) - G(y_1 + x))\, dx + \left( \frac{\lambda}{\mu^2} + \lambda\beta^2 \right)(s_2 - s_1) \\
& \quad = \frac{\lambda}{\mu^2} \left( \int_{y_2}^{y_2 + s_1} - \int_{y_1}^{y_1 + s_1} \right) G(x)\, dx + \left( \frac{\lambda}{\mu^2} + \lambda\beta^2 \right)(s_2 - s_1) \\
& \quad \leq \frac{\lambda}{\mu^2} |y_2 - y_1| + \left( \frac{\lambda}{\mu^2} + \lambda\beta^2 \right)(s_2 - s_1).
\end{aligned}
$$

Similarly, if $y_1 > y_2$, then the right-hand side of (5.37) is equal to

$$\lambda \int_0^{s_1} \left( \frac{1}{\mu^2} G(y_2 + x) + \beta^2 \right) (G(y_1 + x) - G(y_2 + x)) \, dx$$

$$+ \lambda \int_{s_1}^{s_2} \left( \frac{1}{\mu^2} G(y_2 + x) + \beta^2 \right) (1 - G(y_2 + x)) \, dx$$

$$\leq \lambda \left( \frac{1}{\mu^2} + \beta^2 \right) \int_0^{s_1} (G(y_1 + x) - G(y_2 + x)) \, dx + \lambda \left( \frac{1}{\mu^2} + \beta^2 \right) (s_2 - s_1)$$

$$\leq \lambda \left( \frac{1}{\mu^2} + \beta^2 \right) (|y_2 - y_1| + (s_2 - s_1)).$$

Thus,

$$(5.38) \quad |\mathbb{E} Y(s_2, y_2)^2 - \mathbb{E}[Y(s_1, y_1) Y(s_2, y_2)]| \leq C_1 \|(s_1, y_1) - (s_2, y_2)\|,$$

where $C_1 = \sqrt{2} \lambda (\frac{1}{\mu^2} + \beta^2)$. By a similar (in fact, simpler) argument, we get

$$(5.39) \quad |\mathbb{E} Y(s_1, y_1)^2 - \mathbb{E}[Y(s_1, y_1) Y(s_2, y_2)]| \leq C_1 \|(s_1, y_1) - (s_2, y_2)\|.$$

Relations (5.36), (5.38) and (5.39) give (5.35) with $\tilde{C} = 2C_1$.  □

LEMMA 5.2.   *Let $\mathcal{F}_n \triangleq \sigma(u_j^{(n)}, j = 1, \ldots)$ and let $K^{(n)}$ be $\mathcal{F}_n$-measurable stochastic processes such that $K^{(n)} \Rightarrow 0$ in $D[0, \infty)$. For $s \geq 0$ and $y \leq y^*$, let*

$$Z^{(n)}(s, y) \triangleq \frac{1}{n^{1/4}} \sum_{j=1}^{\lfloor \lambda^{(n)} s \sqrt{n} \rfloor} \left( M_j^{(n)}(y^*) - M_j^{(n)} \left( y + \frac{j}{\lambda^{(n)} \sqrt{n}} - K^{(n)} \left( \frac{j}{\lambda^{(n)} \sqrt{n}} \right) \right) \right).$$

*Then, for every $s_0 > 0$ and $y_0 < y^*$, $Z^{(n)} - Y^{(n)} \Rightarrow 0$ in $D([0, s_0] \times [y_0, y^*])$.*

PROOF.   Because $K^{(n)} \Rightarrow 0$, there exists a sequence $\varepsilon_n$ of positive numbers converging to zero so slowly that $\mathbb{P}(A_n) \to 1$, where $A_n = [\sup_{0 \leq s \leq s_0} |K^{(n)}(s)| \leq \varepsilon_n]$. It suffices to prove that

$$(5.40) \qquad\qquad (Z^{(n)} - Y^{(n)}) \mathbb{I}_{A_n} \Rightarrow 0$$

in $D([0, s_0] \times [y_0, y^*])$. Relation (5.40) is a statement about *weak convergence* of stochastic processes, so the underlying probability spaces are irrelevant. Thus, without loss of generality, we can assume that all the random variables (and, thus, all the prelimit processes) under consideration are defined on the same probability space $(\Omega, \mathcal{A}, \mathbb{P})$ and, moreover, *all* the arrival times $\{u_j^{(n)}\}_{j,n=1}^\infty$ are independent of *all* the service times $\{v_j^{(n)}\}_{j,n=1}^\infty$ and *all* the lead times $\{L_j^{(n)}\}_{j,n=1}^\infty$. This is not a limiting assumption, because if, for different $n$, the probability spaces $(\Omega^{(n)}, \mathcal{A}^{(n)}, \mathbb{P}^{(n)})$ on which the sequences $\{u_j^{(n)}\}_{j=1}^{\infty}$, $\{v_j^{(n)}\}_{j=1}^{\infty}$, $\{L_j^{(n)}\}_{j=1}^{\infty}$ are defined are different, we can take $(\Omega, \mathcal{A}, \mathbb{P}) = \prod_{n=1}^{\infty}(\Omega^{(n)}, \mathcal{A}^{(n)}, \mathbb{P}^{(n)})$. Let $\mathcal{F} \triangleq \sigma(u_j^{(n)}, j, n = 1, \ldots)$. Note that, for any (deterministic) function $f$,

$$(5.41) \qquad \mathbb{E}[f(Y^{(n)}, Z^{(n)})|\mathcal{F}] = \mathbb{E}[f(Y^{(n)}, Z^{(n)})|K^{(n)}(\cdot)],$$

that is, we can evaluate the left-hand side of (5.41) by conditioning on a sample path of $K^{(n)}(\cdot)$. As in the proof of Proposition 5.1, let $U^{(n)}$, $n = 1, \ldots$, be a sequence of random variables with distribution (5.17), let $X$ be a random variable independent of this sequence and having cumulative distribution function $G$, and let $(X, U^{(1)}, U^{(2)}, \ldots)$ be independent of $\mathcal{F}$. Let $\overline{m}_2^{(n)}$ be a random measure, depending on the sample path of $K^{(n)}(\cdot)$, defined by

$$
\begin{aligned}
&\overline{m}_2^{(n)}((s_1, s_2] \times (y_1, y_2]) \\
&\triangleq \mathbb{P}[(U^{(n)}, X - U^{(n)} + K^{(n)}(U^{(n)})) \in (s_1, s_2] \times (y_1, y_2]|K^{(n)}(\cdot)] \\
&= \frac{1}{\lfloor \lambda^{(n)}\sqrt{n}s_0 \rfloor} \sum_{j=\lfloor s_1\lambda^{(n)}\sqrt{n} \rfloor + 1}^{\lfloor s_2\lambda^{(n)}\sqrt{n} \rfloor} (G_j^{(n)}(y_2) - G_j^{(n)}(y_1)),
\end{aligned}
$$

where $0 \leq s_1 < s_2 \leq s_0$, $y_0 \leq y_1 < y_2 \leq y^*$ and

$$G_j^{(n)}(y) \triangleq G\left(y + \frac{j}{\lambda^{(n)}\sqrt{n}} - K^{(n)}\left(\frac{j}{\lambda^{(n)}\sqrt{n}}\right)\right).$$

Let us also define random measures $m_2^{(n)} \triangleq \overline{m}_2^{(n)} \mathbb{I}_{A_n} + m_1 \mathbb{I}_{A_n^c}$, where $m_1$ is the measure defined by (5.25). Observe that, for every $\omega \in \Omega$,

$$(5.42) \qquad m_2^{(n)}(\omega) \Rightarrow m_1.$$

Proceeding as in the proof of Proposition 5.1, we get

$$
\begin{aligned}
(5.43) \qquad &\mathbb{E}[(Z^{(n)}(B_1))^2 (Z^{(n)}(B_2))^2 \mathbb{I}_{A_n}|K^{(n)}(\cdot)] \\
&\leq C_1^2 \overline{m}_2^{(n)}(B_1)\overline{m}_2^{(n)}(B_2)\mathbb{I}_{A_n} \\
&\leq C_1^2 m_2^{(n)}(B_1)m_2^{(n)}(B_2)
\end{aligned}
$$

for some deterministic constant $C_1$ and any "neighboring blocks" $B_1, B_2 \subseteq [0, s_0] \times [y_0, y^*]$. As in the proof of Proposition 5.1, (5.42)–(5.43) imply that the random fields $Z^{(n)}\mathbb{I}_{A_n}$ are conditionally tight with respect to $K^{(n)}(\cdot)$. By Proposition 5.1 and the fact that $Y^{(n)}$ is independent of $\mathcal{F}$, the random fields $(Z^{(n)} - Y^{(n)})\mathbb{I}_{A_n}$ are also conditionally tight with respect to $K^{(n)}(\cdot)$. For any $0 \leq s \leq s_0$ and $y_0 \leq y \leq y^*$, we have

$$\mathbb{E}[(Z^{(n)} - Y^{(n)})^2(s, y)\mathbb{I}_{A_n}|K^{(n)}(\cdot)]$$

$$(5.44) \quad \leq \frac{\mathbb{E}(v_1^{(n)})^2}{\sqrt{n}} \sum_{j=1}^{\lfloor \lambda^{(n)} s\sqrt{n} \rfloor} \left( G\left(y + \frac{j}{\lambda^{(n)}\sqrt{n}} + \varepsilon_n\right) \right.$$

$$\left. - G\left(y + \frac{j}{\lambda^{(n)}\sqrt{n}} - \varepsilon_n\right) \right)$$

$$\leq \tilde{C} m_1^{(n)}((0, s] \times (y - \varepsilon_n, y + \varepsilon_n]) \to 0,$$

where $\tilde{C}$ is a constant independent of $n$, because, as we have seen in the proof of Proposition 5.1, the measures $m_1^{(n)}$ converge weakly to a continuous measure $m_1$. In particular, the finite-dimensional conditional distributions of $(Z^{(n)} - Y^{(n)})\mathbb{I}_{A_n}$ with respect to $K^{(n)}(\cdot)$ converge to zero. Thus, the conditional distributions of the random fields $(Z^{(n)} - Y^{(n)})\mathbb{I}_{A_n}$ on $D([0, s_0] \times [y_0, y^*])$ with respect to $K^{(n)}(\cdot)$ converge weakly to 0, so for any continuous and bounded function $f: D([0, s_0] \times [y_0, y^*]) \to \mathbb{R}$,

$$\mathbb{E}[f((Z^{(n)} - Y^{(n)})\mathbb{I}_{A_n})|\mathcal{F}] = \mathbb{E}[f((Z^{(n)} - Y^{(n)})\mathbb{I}_{A_n})|K^{(n)}(\cdot)] \to f(0),$$

where the equality follows from (5.41). This implies (5.40), because, by the bounded convergence theorem,

$$\mathbb{E}f((Z^{(n)} - Y^{(n)})\mathbb{I}_{A_n}) = \mathbb{E}\{\mathbb{E}[f((Z^{(n)} - Y^{(n)})\mathbb{I}_{A_n})|\mathcal{F}]\} \to f(0). \qquad \square$$

PROPOSITION 5.3. $R^{(n)} \Rightarrow 0$ in $D(-\infty, y^*]$.

PROOF. Let $y \leq y^*$ be given. For any $1 \leq j \leq A^{(n)}(nt) - A^{(n)}(nt - \sqrt{n} \times (y^* - y))$ and $y \leq \ell \leq y^*$, $A^{(n)}(nt) - A^{(n)}(nt - \sqrt{n}(\ell - y)) \geq j$ if and only if $nt - \sqrt{n}(\ell - y) < S^{(n)}_{A^{(n)}(nt)-j+1}$ which, in turn, is equivalent to $\ell > y + \frac{1}{\sqrt{n}}(nt - S^{(n)}_{A^{(n)}(nt)-j+1})$. Thus,

$$\int_y^{y^*} \sum_{j=1}^{A^{(n)}(nt) - A^{(n)}(nt - \sqrt{n}(\ell - y))} dM_j^{(n)}(\ell)$$

$$(5.45) \quad = \sum_{j=1}^{A^{(n)}(nt) - A^{(n)}(nt - \sqrt{n}(y^* - y))} \left( M_j^{(n)}(y^*) \right.$$

$$\left. - M_j^{(n)}\left(y + \frac{1}{\sqrt{n}}(nt - S^{(n)}_{A^{(n)}(nt)-j+1})\right) \right).$$

For $s \geq 0$ and $y \leq y^*$, let

$$X^{(n)}(s, y) \triangleq \frac{1}{n^{1/4}} \sum_{j=1}^{\lfloor \lambda^{(n)} s\sqrt{n} \rfloor} \left( M_j^{(n)}(y^*) - M_j^{(n)}\left(y + \frac{1}{\sqrt{n}}(nt - S^{(n)}_{(A^{(n)}(nt)-j+1)^+})\right) \right)$$

$$= \frac{1}{n^{1/4}} \sum_{j=1}^{\lfloor \lambda^{(n)} s \sqrt{n} \rfloor} \left( M_j^{(n)}(y^*) - M_j^{(n)}\left( y + \frac{j}{\lambda^{(n)}\sqrt{n}} - H^{(n)}\left( \frac{j}{\lambda^{(n)}\sqrt{n}} \right) \right) \right),$$

where, for $u \geq 0$,

$$H^{(n)}(u) \triangleq \frac{1}{\sqrt{n}} [S^{(n)}_{\lfloor A^{(n)}(nt) + 1 - \lambda^{(n)}\sqrt{n}u \rfloor \vee 0} - (nt - \sqrt{n}u)]$$

$$= \widehat{S}^{(n)}_{((1/n)A^{(n)}(nt) + (1/n) - \lambda^{(n)}u/\sqrt{n})^+} + \frac{1}{\lambda^{(n)}\sqrt{n}} \lfloor A^{(n)}(nt) + 1 - \lambda^{(n)}\sqrt{n}u \rfloor^+$$

$$\quad - \sqrt{n}t + u.$$

Observe that, by (2.15), (2.16) and the differencing theorem, with probability approaching 1,

$$H^{(n)}(u) = \widehat{S}^{(n)}_{\lambda^{(n)}t + o(1)} + \frac{1}{\lambda^{(n)}} \widehat{A}^{(n)}(t) + O\left( \frac{1}{\sqrt{n}} \right) \Rightarrow 0.$$

By Lemma 5.2, for every $s_0 > 0$ and $y_0 < y^*$,

(5.46)                          $X^{(n)} - Y^{(n)} \Rightarrow 0$

in $D([0, s_0] \times [y_0, y^*])$ and, thus, by Proposition 5.1,

(5.47)                                $X^{(n)} \Rightarrow Y$

in $D([0, s_0] \times [y_0, y^*], \mathbb{R})$. By (5.10), (5.45), the definition of $X^{(n)}$, (5.6) and (5.9),

$$R^{(n)}(y) = X^{(n)}\left( y^* - y + \frac{1}{\lambda^{(n)}n^{1/4}} \widetilde{A}^{(n)}(y^* - y), y \right) - Y^{(n)}(y^* - y, y-)$$

$$= \left( X^{(n)}\left( y^* - y + \frac{1}{\lambda^{(n)}n^{1/4}} \widetilde{A}^{(n)}(y^* - y), y \right) - X^{(n)}(y^* - y, y) \right)$$

$$\quad + (X^{(n)}(y^* - y, y) - Y^{(n)}(y^* - y, y))$$

$$\quad + (Y^{(n)}(y^* - y, y) - Y^{(n)}(y^* - y, y-)),$$

which converges to 0 in $D[y_0, y^*]$ for every $y_0 < y^*$ by (4.11), Proposition 5.1, (5.46), (5.47), continuity of the sample paths of $Y$ and the differencing theorem.  $\square$

5.2. *Asymptotic analysis for* $I_2^{(n)}$.  The analysis of the limiting behavior of $n^{1/4} I_2^{(n)}(y)$ is easier. We have the following:

LEMMA 5.4.  *We have* $n^{1/4}(I_2^{(n)} - H) \Rightarrow Z$ *in* $D(-\infty, y^*]$, *where* $Z$ *is a mean-zero Gaussian process with continuous paths and covariance*

$$\mathbb{E}[Z(y_1)Z(y_2)] = \lambda \alpha^2 \rho^2 \int_{y_1}^{y^*} \int_{y_2}^{y^*} (\ell - y_1) \wedge (k - y_2) \, dG(\ell) \, dG(k).$$

PROOF.    For $y \le y^*$, let $Z(y) \triangleq \frac{1}{\mu} \int_y^{y^*} B(\ell - y) \, dG(\ell)$, where $B$ is a zero-drift Brownian motion with variance $\alpha^2 \lambda^3$ per unit time. It is easy to see that $Z$ is a mean-zero Gaussian process with continuous sample paths. Furthermore,

$$\mathbb{E}[Z(y_1)Z(y_2)] = \mathbb{E}\left[ \frac{1}{\mu} \int_{y_1}^{y^*} B(\ell - y_1) \, dG(\ell) \cdot \frac{1}{\mu} \int_{y_2}^{y^*} B(k - y_2) \, dG(k) \right]$$

$$= \frac{1}{\mu^2} \int_{y_1}^{y^*} \int_{y_2}^{y^*} \mathbb{E}[B(\ell - y_1)B(k - y_2)] \, dG(\ell) \, dG(k)$$

$$= \alpha^2 \lambda \rho^2 \int_{y_1}^{y^*} \int_{y_2}^{y^*} (\ell - y_1) \wedge (k - y_2) \, dG(\ell) \, dG(k).$$

Finally, by (5.5),

$$I_2^{(n)}(y) = \frac{1}{\mu^{(n)}\sqrt{n}} \int_y^{y^*} (A^{(n)}(nt) - A^{(n)}(nt - \sqrt{n}(\ell - y))) \, dG(\ell)$$

$$(5.48) \qquad = \frac{1}{\mu^{(n)}} \int_y^{y^*} \left( \lambda^{(n)}(\ell - y) + \widehat{A}^{(n)}(t) - \widehat{A}^{(n)}\left( t - \frac{\ell - y}{\sqrt{n}} \right) \right) dG(\ell)$$

$$= \rho^{(n)} \int_y^{y^*} (\ell - y) \, dG(\ell) + \int_y^{y^*} \frac{1}{\mu^{(n)} n^{1/4}} \widetilde{A}^{(n)}(\ell - y) \, dG(\ell),$$

and thus, by (2.9), (2.10), (4.11) and the fact that

$$(5.49) \qquad\qquad \int_y^{y^*} (\ell - y) \, dG(\ell) = H(y),$$

we have [with $\widetilde{A}$ defined by (4.11)]

$$(5.50) \qquad n^{1/4}(I_2^{(n)}(y) - H(y))$$

$$= \int_y^{y^*} \frac{1}{\mu^{(n)}} \widetilde{A}^{(n)}(\ell - y) \, dG(\ell) + O\left( \frac{1}{n^{1/4}} \right) \Rightarrow Z(y). \qquad \square$$

### 5.3. Asymptotic analysis for $(\widehat{J}^{(n)}, \widehat{W}^{(n)}(t))$.

COROLLARY 5.5.    We have $\widehat{J}^{(n)} \Rightarrow J^*$ in $D(-\infty, y^*]$, where $J^*$ is a mean-zero Gaussian process with continuous paths and covariance (3.5).

PROOF.    By (5.7), (5.9), (5.10), (5.50), Propositions 5.1, 5.3 and the independence of the arrivals, the service times and the lead times,

$$(n^{1/4}I_1^{(n)}(y), n^{1/4}(I_2^{(n)}(y) - H(y)))$$

$$= \left( \frac{1}{n^{1/4}} \int_y^{y^* } \sum_{j=0}^{\lambda^{(n)}\sqrt{n}(\ell - y) + n^{1/4}\widetilde{A}^{(n)}(\ell - y) - 1} dM_{A^{(n)}(nt) - j}^{(n)}(\ell), \right.$$

$$\left. \int_y^{y^*} \frac{1}{\mu^{(n)}} \widetilde{A}^{(n)}(\ell - y)\, dG(\ell) \right) + o(1)$$

$$(5.51) \quad \stackrel{d}{=} \left( \frac{1}{n^{1/4}} \int_y^{y^*} \sum_{j=1}^{\lambda^{(n)}\sqrt{n}(\ell - y) + n^{1/4}\widetilde{A}^{(n)}(\ell - y)} dM_j^{(n)}(\ell), \right.$$

$$\left. \int_y^{y^*} \frac{1}{\mu^{(n)}} \widetilde{A}^{(n)}(\ell - y)\, dG(\ell) \right) + o(1)$$

$$= \left( Y^{(n)}(y^* - y, y-) + R^{(n)}(y), \int_y^{y^*} \frac{1}{\mu^{(n)}} \widetilde{A}^{(n)}(\ell - y)\, dG(\ell) \right) + o(1)$$

$$\Rightarrow (Y(y^* - y, y), Z(y)),$$

where $Y$ and $Z$ are as in Proposition 5.1 and Lemma 5.4, independent of each other. Thus, by (5.2), (5.3) and (5.51), $\widehat{J}^{(n)} \Rightarrow J^*$, as claimed.  $\square$

PROPOSITION 5.6.   *We have* $(\widehat{J}^{(n)}, \widehat{W}^{(n)}(t)) \Rightarrow (J^*, W^*(t))$ *in* $D(-\infty, y^*] \times \mathbb{R}$, *where* $J^*$ *is as in Corollary* 5.5 *and* $W^*$ *is a reflected Brownian motion with variance* $(\alpha^2 \rho^2 + \beta^2)\lambda$ *per unit time and drift* $-\gamma$, *independent of* $J^*$.

PROOF.   Fix an arbitrary $y_0 < y^*$. By (2.17), Corollary 5.5 and its proof, it suffices to show that $\widehat{W}^{(n)}(t)$ is asymptotically independent of the pair of processes $(n^{1/4} I_1^{(n)}(y), n^{1/4}(I_2^{(n)}(y) - H(y)))$, $y_0 \le y \le y^*$. We assume throughout the proof that $n$ is sufficiently large so that $nt - \sqrt{n}(y^* - y_0) > 0$. We note at the outset that, by (2.17) and the differencing theorem,

$$(5.52) \qquad \widehat{W}^{(n)}\left( t - \frac{1}{\sqrt{n}}(y^* - y_0) \right) - \widehat{W}^{(n)}(t) \Rightarrow 0.$$

Denote $T_1^{(n)} \triangleq S_{A^{(n)}(nt - \sqrt{n}(y^* - y_0))}^{(n)}$ and $T_2^{(n)} \triangleq S_{A^{(n)}(nt - \sqrt{n}(y^* - y_0)) + 1}^{(n)}$, so that $T_1^{(n)} \le nt - \sqrt{n}(y^* - y_0) < T_2^{(n)}$. Define $\theta^{(n)} = T_2^{(n)} - (nt - \sqrt{n}(y^* - y_0))$. The process $A^{(n)}(s)$, $0 \le s \le nt - \sqrt{n}(y^* - y_0)$, and the process

$$\overline{A}^{(n)}(y) \triangleq A^{(n)}(T_2^{(n)} + \sqrt{n}(y - y_0)) - A^{(n)}(T_2^{(n)})$$
$$= A^{(n)}(nt - \sqrt{n}(y^* - y) + \theta^{(n)}) - A^{(n)}(nt - \sqrt{n}(y^* - y_0) + \theta^{(n)}),$$

$y_0 \le y \le y^*$, are independent.

We show now that

$$(5.53) \qquad\qquad \frac{1}{n^{1/4}} \theta^{(n)} \Rightarrow 0.$$

Using the process $\widehat{S}^{(n)}$ that has the continuous limit $S^*$ in (2.15), we may write

$$\max_{0 \leq \tau \leq 2\lambda t} |\widehat{S}^{(n)}(\tau) - \widehat{S}^{(n)}(\tau-)| = \frac{1}{\sqrt{n}} \max_{1 \leq j \leq \lfloor 2\lambda nt \rfloor} \left| u_j^{(n)} - \frac{1}{\lambda^{(n)}} \right|$$

$$\geq \max_{1 \leq j \leq \lfloor 2\lambda nt \rfloor} \frac{1}{\sqrt{n}} u_j^{(n)} - \frac{1}{\sqrt{n}\lambda^{(n)}}.$$

Because the limit of $\widehat{S}^{(n)}$ is continuous, the term $\max_{1 \leq j \leq \lfloor 2\lambda nt \rfloor} \frac{1}{\sqrt{n}} u_j^{(n)}$ converges to zero in probability. Since $\widehat{A}^{(n)} \Rightarrow A^*$, we can choose a sequence of sets $\{B_n\}_{n=1}^{\infty}$ with $\mathbb{P}(B_n) \to 1$ such that $A^{(n)}(nt) \leq 2\lambda nt - 1$ on $B_n$. We set

$$s_1^{(n)} = \frac{1}{\sqrt{n}}(nt - T_1^{(n)}), \qquad s_2^{(n)} = \frac{1}{\sqrt{n}}(nt - T_2^{(n)}),$$

so

$$\frac{1}{\sqrt{n}} u_{A^{(n)}(nt - \sqrt{n}(y^* - y_0)) + 1} = \frac{1}{\sqrt{n}}(T_2^{(n)} - T_1^{(n)}) = s_1^{(n)} - s_2^{(n)} \geq \frac{1}{\sqrt{n}} \theta^{(n)}.$$

In particular, $(s_1^{(n)} - s_2^{(n)}) \mathbb{I}_{B_n}$ converges to zero in probability, and hence, so does $s_1^{(n)} - s_2^{(n)}$. We have $s_2^{(n)} = s_1^{(n)} - (s_1^{(n)} - s_2^{(n)}) \geq y^* - y_0 - (s_1^{(n)} - s_2^{(n)})$, so $\mathbb{P}(C_n) \to 1$, where $C_n = \lceil s_2^{(n)} \geq 0 \rceil$. On the set $C_n$, the differencing theorem now implies that zero is the limit in probability of

$$\widetilde{A}^{(n)}(s_1^{(n)}) - \widetilde{A}^{(n)}(s_2^{(n)}) = \frac{1}{n^{1/4}}(A^{(n)}(nt - s_2^{(n)}\sqrt{n}) - A^{(n)}(nt - s_1^{(n)}\sqrt{n})$$

$$- \lambda^{(n)}\sqrt{n}(s_1^{(n)} - s_2^{(n)}))$$

$$= \frac{1}{n^{1/4}}(1 - \lambda^{(n)}\sqrt{n}(s_1^{(n)} - s_2^{(n)})).$$

We conclude that $n^{1/4}(s_1^{(n)} - s_2^{(n)}) \Rightarrow 0$. This implies (5.53).

Recall that

$$n^{1/4} I_1^{(n)}(y) = \frac{1}{n^{1/4}} \int_y^{y^*} \sum_{j=A^{(n)}(nt - \sqrt{n}(\ell - y)) + 1}^{A^{(n)}(nt)} dM_j^{(n)}(\ell), \qquad y_0 \leq y \leq y^*.$$

We define the related process $\widetilde{I}_1^{(n)}(y)$ by

$$n^{1/4} \widetilde{I}_1^{(n)}(y) = \frac{1}{n^{1/4}} \int_y^{y^*} \sum_{j=A^{(n)}(nt - \sqrt{n}(\ell - y) + \theta^{(n)}) + 1}^{A^{(n)}(nt + \theta^{(n)})} dM_j^{(n)}(\ell)$$

$$= \frac{1}{n^{1/4}} \int_y^{(y + \theta^{(n)}/\sqrt{n}) \wedge y^*} \sum_{j=A^{(n)}(nt - \sqrt{n}(\ell - y) + \theta^{(n)}) + 1}^{A^{(n)}(nt + \theta^{(n)})} dM_j^{(n)}(\ell)$$

(5.54)

$$+ \frac{1}{n^{1/4}} \int_{(y+\theta^{(n)}/\sqrt{n}) \wedge y^*}^{y^*} \sum_{j=A^{(n)}(nt)+1}^{A^{(n)}(nt+\theta^{(n)})} dM_j^{(n)}(\ell)$$

$$+ n^{1/4} I_1^{(n)} \left( \left( y + \frac{\theta^{(n)}}{\sqrt{n}} \right) \wedge y^* \right).$$

Recall further that, by (5.48) and (5.49),

$$n^{1/4}(I_2^{(n)}(y) - H(y)) = \frac{n^{1/4}}{\mu^{(n)}} \int_y^{y^*} \left[ \widehat{A}^{(n)}(t) - \widehat{A}^{(n)} \left( t - \frac{\ell-y}{\sqrt{n}} \right) \right] dG(\ell)$$

$$+ n^{1/4}(\rho^{(n)} - 1) H(y), \qquad y_0 \le y \le y^*.$$

We define the related process $\tilde{I}_2^{(n)}(y)$ by

$$n^{1/4}(\tilde{I}_2^{(n)}(y) - H(y))$$

$$= \frac{n^{1/4}}{\mu^{(n)}} \int_y^{y^*} \left[ \widehat{A}^{(n)} \left( t + \frac{\theta^{(n)}}{n} \right) - \widehat{A}^{(n)} \left( t - \frac{\ell-y}{\sqrt{n}} + \frac{\theta^{(n)}}{n} \right) \right] dG(\ell)$$

$$= \frac{n^{1/4}}{\mu^{(n)}} \int_y^{(y+\theta^{(n)}/\sqrt{n}) \wedge y^*} \left[ \widehat{A}^{(n)} \left( t + \frac{\theta^{(n)}}{n} \right) \right.$$

(5.55)

$$\left. - \widehat{A}^{(n)} \left( t - \frac{\ell-y}{\sqrt{n}} + \frac{\theta^{(n)}}{n} \right) \right] dG(\ell)$$

$$+ n^{1/4} \left( I_2^{(n)} \left( \left( y + \frac{\theta^{(n)}}{\sqrt{n}} \right) \wedge y^* \right) - H \left( y + \frac{\theta^{(n)}}{\sqrt{n}} \right) \right)$$

$$- n^{1/4}(\rho^{(n)} - 1) H \left( y + \frac{\theta^{(n)}}{\sqrt{n}} \right)$$

$$+ \frac{n^{1/4}}{\mu^{(n)}} \left[ \widehat{A}^{(n)} \left( t + \frac{\theta^{(n)}}{n} \right) - \widehat{A}^{(n)}(t) \right] \left[ G(y^*) - G \left( y + \frac{\theta^{(n)}}{\sqrt{n}} \right) \right].$$

Note in (5.55) that $G(y) = G(y^*) = 1$ and $H(y) = H(y^*) = 0$ for $y \ge y^*$, so evaluating $G$ and $H$ at $y + \frac{\theta^{(n)}}{\sqrt{n}}$ gives the same result as evaluating these functions at $(y + \frac{\theta^{(n)}}{\sqrt{n}}) \wedge y^*$. The pair of processes $(\tilde{I}_1^{(n)}(y), \tilde{I}_2^{(n)}(y); y_0 \le y \le y^*)$ is independent of the random variable $\widehat{W}^{(n)}(t - \frac{1}{\sqrt{n}}(y^* - y_0))$, as we now explain. In $\tilde{I}_1^{(n)}(y)$, the sum $j = A^{(n)}(nt - \sqrt{n}(\ell - y) + \theta^{(n)}) + 1$ to $A^{(n)}(nt + \theta^{(n)})$ has

$$A^{(n)}(nt + \theta^{(n)}) - A^{(n)}(nt - \sqrt{n}(\ell - y) + \theta^{(n)}) = \overline{A}^{(n)}(y^*) - \overline{A}^{(n)}(y^* - (l - y))$$

terms, and this is independent of $\widehat{W}^{(n)}(t - \frac{1}{\sqrt{n}}(y^* - y_0))$. The $M_j^{(n)}$ processes appearing in this sum involve service times $v_j^{(n)}$, but these particular indices $j$ do not appear in the definition of $\widehat{W}^{(n)}(t - \frac{1}{\sqrt{n}}(y^* - y_0))$. The integrand appearing in $\tilde{I}_2^{(n)}(y)$,

$$\widehat{A}^{(n)}\left(t + \frac{\theta^{(n)}}{n}\right) - \widehat{A}^{(n)}\left(t - \frac{\ell - y}{\sqrt{n}} + \frac{\theta^{(n)}}{n}\right)$$

$$= \frac{1}{\sqrt{n}}[A^{(n)}(nt + \theta^{(n)}) - A^{(n)}(nt - \sqrt{n}(\ell - y) + \theta^{(n)}) - \lambda^{(n)}\sqrt{n}(l - y)]$$

$$= \frac{1}{\sqrt{n}}[\overline{A}^{(n)}(y^*) - \overline{A}^{(n)}(y^* - (l - y)) - \lambda^{(n)}\sqrt{n}(l - y)],$$

is also independent of $\widehat{W}^{(n)}(t - \frac{1}{\sqrt{n}}(y^* - y_0))$.

It remains to show

$$(5.56) \quad n^{1/4}(\tilde{I}_1^{(n)}(y) - I_1^{(n)}(y)) \Rightarrow 0, \qquad n^{1/4}(\tilde{I}_2^{(n)}(y) - I_2^{(n)}(y)) \Rightarrow 0.$$

This will imply that the limit of the pair $(n^{1/4}I_1^{(n)}(y), n^{1/4}(I_2^{(n)}(y) - H(y)))$ is independent of $W^*(t)$, the limit of $\widehat{W}^{(n)}(t - \frac{1}{\sqrt{n}}(y^* - y_0))$. Since $\widehat{J}^{(n)}(y) = n^{1/4}I_1^{(n)}(y) + n^{1/4}(I_2^{(n)}(y) - H(y)) \Rightarrow J^*(y)$, we will have the desired result.

From (5.54), we have

$$
\begin{aligned}
& n^{1/4}(\tilde{I}_1^{(n)}(y) - I_1^{(n)}(y)) \\
&= \frac{1}{n^{1/4}} \int_y^{(y + \theta^{(n)}/\sqrt{n}) \wedge y^*} \sum_{j=A^{(n)}(nt - \sqrt{n}(\ell - y) + \theta^{(n)}) + 1}^{A^{(n)}(nt + \theta^{(n)})} dM_j^{(n)}(\ell) \\
&\quad + \frac{1}{n^{1/4}} \sum_{j=A^{(n)}(nt)+1}^{A^{(n)}(nt+\theta^{(n)})} \left[ M_j^{(n)}(y^*) - M_j^{(n)}\left( \left(y + \frac{\theta^{(n)}}{\sqrt{n}}\right) \wedge y^* \right) \right] \\
&\quad + n^{1/4}\left[ I_1^{(n)}\left( \left(y + \frac{\theta^{(n)}}{\sqrt{n}}\right) \wedge y^* \right) - I_1^{(n)}(y) \right].
\end{aligned}
$$

(5.57)

For every $y \in [y_0, y^*]$, the absolute value of each of the first two terms on the right-hand side of (5.57) is bounded above by

$$(5.58) \qquad \frac{1}{n^{1/4}} \sum_{j=A^{(n)}(nt)+1}^{A^{(n)}(nt+\theta^{(n)})} \left( v_j^{(n)} + \frac{1}{\mu^{(n)}} \right).$$

The ordinary and renewal functional central limit theorems for triangular arrays (see, e.g., [5, 7, 13] and Theorem 14.6 in [3]) imply that

$$
\begin{aligned}
(5.59) \quad \widetilde{C}^{(n)}(s) &\triangleq \frac{1}{n^{1/4}}[A^{(n)}(nt + \sqrt{n}s) - A^{(n)}(nt) - \lambda^{(n)}\sqrt{n}s] \\
&= n^{1/4}\left[\widehat{A}^{(n)}\left(t + \frac{s}{\sqrt{n}}\right) - \widehat{A}^{(n)}(t)\right] \Rightarrow C^*(s),
\end{aligned}
$$

where $C^*$ is a Brownian motion with zero drift and variance $\alpha^2\lambda^3$ per unit time. In particular, $\widetilde{C}^{(n)}(sn^{-1/4}) \Rightarrow 0$ in $D[0,\infty)$ and, thus,

$$
(5.60) \quad \frac{1}{n^{1/4}}[A^{(n)}(nt + sn^{1/4}) - A^{(n)}(nt)] \Rightarrow \lambda s
$$

in $D[0,\infty)$. Let $\varepsilon > 0$ be given. The convergence (5.53) implies that $\mathbb{P}(D_n) \to 1$, where $D_n = [\theta^{(n)} \leq \varepsilon n^{1/4}]$. On $D_n$, the expression (5.58) is dominated by

$$
\begin{aligned}
&\frac{1}{n^{1/4}} \sum_{j=A^{(n)}(nt)+1}^{A^{(n)}(nt+\varepsilon n^{1/4})} \left(v_j^{(n)} + \frac{1}{\mu^{(n)}}\right) \\
&= \frac{A^{(n)}(nt + \varepsilon n^{1/4}) - A^{(n)}(nt)}{n^{1/4}} \\
&\quad \times \frac{1}{A^{(n)}(nt + \varepsilon n^{1/4}) - A^{(n)}(nt)} \sum_{j=A^{(n)}(nt)+1}^{A^{(n)}(nt+\varepsilon n^{1/4})} \left(v_j^{(n)} + \frac{1}{\mu^{(n)}}\right),
\end{aligned}
$$

which converges weakly to $\lambda\varepsilon \cdot \frac{2}{\lambda} = 2\varepsilon$ by (5.60) and the law of large numbers for triangular arrays, together with the independence of the arrival times and the service times. Thus, the expression (5.58) converges to zero in probability. To see that the third term on the right-hand side of (5.57) converges to zero, we use (5.53), the fact that $n^{1/4}I_1^{(n)}(y) \Rightarrow Y(y^* - y, y)$ by (5.11) and Propositions 5.1, 5.3, the joint continuity of $Y$ and the differencing theorem. This concludes the proof of the first convergence claimed in (5.56).

For the second convergence claimed in (5.56), we use (5.55) to write

$$
\begin{aligned}
(5.61) \quad &n^{1/4}(\widetilde{I}_2^{(n)}(y) - I_2^{(n)}(y)) \\
&= \frac{n^{1/4}}{\mu^{(n)}} \int_y^{(y+\theta^{(n)}/\sqrt{n})\wedge y^*} \left[\widehat{A}^{(n)}\left(t + \frac{\theta^{(n)}}{n}\right)\right. \\
&\qquad\qquad\qquad \left. - \widehat{A}^{(n)}\left(t - \frac{\ell - y}{\sqrt{n}} + \frac{\theta^{(n)}}{n}\right)\right] dG(\ell) \\
&\quad + n^{1/4}\left[I_2^{(n)}\left(\left(y + \frac{\theta^{(n)}}{\sqrt{n}}\right) \wedge y^*\right)\right.
\end{aligned}
$$

$$- H\left(y + \frac{\theta^{(n)}}{\sqrt{n}}\right)\right) - (I_2^{(n)}(y) - H(y))\bigg]$$

$$- n^{1/4}(\rho^{(n)} - 1)H\left(y + \frac{\theta^{(n)}}{\sqrt{n}}\right)$$

$$+ \frac{n^{1/4}}{\mu^{(n)}}\left[\widehat{A}^{(n)}\left(t + \frac{\theta^{(n)}}{n}\right) - \widehat{A}^{(n)}(t)\right]\left[G(y^*) - G\left(y + \frac{\theta^{(n)}}{\sqrt{n}}\right)\right].$$

The absolute value of the first term on the right-hand side of (5.61) is bounded uniformly in $y_0 \leq y \leq y^*$ by $\frac{2}{\mu^{(n)}}\sup_{0 \leq s \leq \theta^{(n)}/\sqrt{n}}|\tilde{C}^{(n)}(s)|$, which converges to zero by (5.53) and (5.59). Convergence to zero of the second term on the right-hand side of (5.61) follows from Lemma 5.4, (5.53) and the differencing theorem. The absolute value of the third term in (5.61) is bounded uniformly in $y_0 \leq y \leq y^*$ by $n^{1/4}|\rho^{(n)} - 1|H(y_0)$, which converges to zero by (2.10). Finally, the last term on the right-hand side of (5.61) converges weakly to zero by (5.53), (5.59) and the differencing theorem. $\square$

COROLLARY 5.7.  *We have* $(\widehat{J}^{(n)}, \widehat{W}^{(n)}(t), \widehat{F}^{(n)}(t)) \Rightarrow (J^*, W^*(t), F^*(t))$ *in the space* $D(-\infty, y^*] \times \mathbb{R}^2$.

PROOF.  By the definition of the frontier, customers with lead times at time $nt$ exceeding $F^{(n)}(nt)$ have not received any service by that time. Thus, by Corollary 3.8 of [4] and the display above (3.32) in [4],

$$(5.62) \quad \begin{aligned} \widehat{W}^{(n)}(t) &= \widehat{\mathcal{V}}^{(n)}(t)(\widehat{F}^{(n)}(t), \infty) + \widehat{\mathcal{W}}^{(n)}(t)[\widehat{C}^{(n)}(t), \widehat{F}^{(n)}(t)] \\ &= \widehat{\mathcal{V}}^{(n)}(t)(\widehat{F}^{(n)}(t), \infty) + o(1) \\ &= H(\widehat{F}^{(n)}(t)) + o(1). \end{aligned}$$

The function $H^{-1}$ is continuous, hence, uniformly continuous on each bounded interval $[0, c]$, $c > 0$, and the Lipschitz constant of $H^{-1}$ on $[c, \infty)$ approaches 1 as $c \to \infty$. In particular, $H^{-1}$ is uniformly continuous on $[0, \infty)$, so (5.62) implies that $H^{-1}(\widehat{W}^{(n)}(t)) - \widehat{F}^{(n)}(t) \Rightarrow 0$. Thus, by Proposition 5.6,

$$\begin{aligned} (\widehat{J}^{(n)}, \widehat{W}^{(n)}(t), \widehat{F}^{(n)}(t)) &= (\widehat{J}^{(n)}, \widehat{W}^{(n)}(t), H^{-1}(\widehat{W}^{(n)}(t))) + o(1) \\ &\Rightarrow (J^*, W^*(t), F^*(t)). \end{aligned} \qquad \square$$

**6. Proofs of the main results.**  In this section we prove Theorems 3.3, 3.4 and Proposition 3.5. We need the following refinement of Lemma 4.1 to include the right endpoint $\widehat{F}^{(n)}(t)$.

LEMMA 6.1.  *We have* $n^{1/4}\widehat{\mathcal{W}}^{(n)}(t)[\widehat{C}^{(n)}(t), \widehat{F}^{(n)}(t)] \Rightarrow 0$.

PROOF. By Proposition 3.1 and Lemma 4.1, it suffices to show that, for any fixed $y_0 < y^*$,

$$(6.1) \qquad n^{1/4} \widehat{\mathcal{W}}^{(n)}(t)\{\widehat{F}^{(n)}(t)\} \mathbb{I}_{\{\widehat{F}^{(n)}(t) \geq y_0\}} \Rightarrow 0.$$

Take a sequence $\varepsilon_n \downarrow 0$ such that $\varepsilon_n = o(n^{-1/4})$. Then

$$
\begin{aligned}
& n^{1/4} \widehat{\mathcal{W}}^{(n)}(t)\{\widehat{F}^{(n)}(t)\} \mathbb{I}_{\{\widehat{F}^{(n)}(t) \geq y_0\}} \\
& \quad \leq n^{1/4} \widehat{\mathcal{V}}^{(n)}(t)\{\widehat{F}^{(n)}(t)\} \mathbb{I}_{\{\widehat{F}^{(n)}(t) \geq y_0\}} \\
& \quad \leq n^{1/4} [\widehat{\mathcal{V}}^{(n)}(t)(\widehat{F}^{(n)}(t) - \varepsilon_n, \infty) - \widehat{\mathcal{V}}^{(n)}(t)(\widehat{F}^{(n)}(t), \infty)] \mathbb{I}_{\{\widehat{F}^{(n)}(t) \geq y_0\}} \\
& \quad \leq [\widehat{J}^{(n)}(\widehat{F}^{(n)}(t) - \varepsilon_n) - \widehat{J}^{(n)}(\widehat{F}^{(n)}(t))] \mathbb{I}_{\{\widehat{F}^{(n)}(t) \geq y_0\}} \\
& \qquad + n^{1/4}(H(\widehat{F}^{(n)}(t) - \varepsilon_n) - H(\widehat{F}^{(n)}(t))) \\
& \quad \leq \omega(\widehat{J}^{(n)}, \varepsilon_n) + n^{1/4} \varepsilon_n,
\end{aligned}
$$

where, for $x \in D(-\infty, y^*]$ and $\delta > 0$,

$$\omega(x, \delta) \triangleq \sup_{\substack{y_0 - \varepsilon_1 \leq s_1 \leq s_2 \leq y^* \\ s_2 - s_1 \leq \delta}} |x(s_2) - x(s_1)|$$

[and we have used Lipschitz continuity of $H$ in the last line of (6.2)]. The right-hand side of the last inequality in (6.2) converges weakly to zero by Corollary 5.5 and the choice of $\varepsilon_n$. This shows (6.1). □

PROOF OF THEOREM 3.3. For every $y \leq y^*$, we have

$$
\begin{aligned}
& n^{1/4}[\widehat{\mathcal{W}}^{(n)}(t)(y, \infty) - H(y \vee \widehat{F}^{(n)}(t))] \\
& \quad = n^{1/4}[\widehat{\mathcal{W}}^{(n)}(t)([\widehat{C}^{(n)}(t), \widehat{F}^{(n)}(t)] \cap (y, \infty)) \\
& \qquad + \widehat{\mathcal{V}}^{(n)}(t)(y \vee \widehat{F}^{(n)}(t), \infty) - H(y \vee \widehat{F}^{(n)})] \\
& \quad = \widehat{J}^{(n)}(y \vee \widehat{F}^{(n)}(t)) + o(1),
\end{aligned}
$$

where the first equality follows from the definition of the frontier as in the proof of Corollary 5.7 and the second one from Lemma 6.1. The mapping $\Phi : D(-\infty, y^*] \times (-\infty, y^*] \rightarrow D(-\infty, y^*]$ defined by $\Phi(x, y)(\cdot) \triangleq x(y \vee \cdot)$ is continuous on $C(-\infty, y^*] \times (-\infty, y^*]$. Thus, by (6.3) and Corollary 5.7,

$$
\begin{aligned}
n^{1/4}[\widehat{\mathcal{W}}^{(n)}(t)(y, \infty) - H(y \vee \widehat{F}^{(n)})] &= \Phi(\widehat{J}^{(n)}, \widehat{F}^{(n)}(t))(y) + o(1) \\
&\Rightarrow \Phi(J^*, F^*(t))(y) \\
&= J^*(y \vee F^*(t))
\end{aligned}
$$

in $D(-\infty, y^*]$. □

PROOF OF THEOREM 3.4. Let $\Psi : [0,\infty)^2 \times \mathbb{R} \to D(-\infty, y^*]$ be the mapping defined by $\Psi(a,a,c)(y) \triangleq c\mathbb{I}_{\{y < H^{-1}(a)\}}$, $\Psi(a,b,c)(y) \triangleq \frac{c}{a-b}(a \wedge H(y) - b \wedge H(y))$ for $a \neq b$. It is easy to check that $\Psi$ is a continuous map of $[0,\infty)^2 \times \mathbb{R}$ into the space $D(-\infty, y^*]$ endowed with the $M_1$ topology. As in (6.3), we have

$$n^{1/4}[\widehat{W}^{(n)}(t) - H(\widehat{F}^{(n)}(t))] = n^{1/4}[\widehat{\mathcal{W}}^{(n)}(t)[\widehat{C}^{(n)}(t), \widehat{F}^{(n)}(t)]$$
$$+ \widehat{\mathcal{Y}}^{(n)}(t)(\widehat{F}^{(n)}(t), \infty) - H(\widehat{F}^{(n)}(t))]$$
$$= \widehat{\mathcal{J}}^{(n)}(\widehat{F}^{(n)}(t)) + o(1),$$

so, by Corollary 5.7, we get

(6.4)
$$(\widehat{W}^{(n)}(t), H(\widehat{F}^{(n)}(t)), n^{1/4}[\widehat{W}^{(n)}(t) - H(\widehat{F}^{(n)}(t))])$$
$$\Rightarrow (W^*(t), W^*(t), J^*(F^*(t))),$$

because the mapping $\Phi : D(-\infty, y^*] \times (-\infty, y^*] \to \mathbb{R}$ defined by $\Phi(x,y) \triangleq x(y)$ is continuous on $C(-\infty, y^*] \times (-\infty, y^*]$. By (1.6) and (6.4),

(6.5)
$$n^{1/4}[H(y \vee H^{-1}(\widehat{W}^{(n)}(t))) - H(y \vee \widehat{F}^{(n)}(t))]$$
$$= n^{1/4}[(H(y) \wedge \widehat{W}^{(n)}(t)) - (H(y) \wedge H(\widehat{F}^{(n)}(t)))]$$
$$= \Psi(\widehat{W}^{(n)}(t), H(\widehat{F}^{(n)}(t)), n^{1/4}[\widehat{W}^{(n)}(t) - H(\widehat{F}^{(n)}(t))])(y)$$
$$\Rightarrow \Psi(W^*(t), W^*(t), J^*(F^*(t)))(y)$$
$$= J^*(F^*(t))\mathbb{I}_{\{y < F^*(t)\}},$$

where weak convergence holds in the $M_1$ topology. The $J_1$ topology on $D(-\infty, y^*]$ is stronger than the $M_1$ topology, so, by Theorem 3.3, (3.4) holds in the $M_1$ topology. It is clear from the above argument and the proof of Theorem 3.3 that the convergence (3.4) and (6.5) is, in fact, joint [because both (3.4) and (6.5) follow from Corollary 5.7 by the continuous mapping theorem]. Thus,

$$n^{1/4}[\widehat{\mathcal{W}}^{(n)}(t)(y, \infty) - H(y \vee H^{-1}(\widehat{W}^{(n)}(t)))]$$
$$= n^{1/4}[\widehat{\mathcal{W}}^{(n)}(t)(y, \infty) - H(y \vee \widehat{F}^{(n)}(t))]$$
$$- n^{1/4}[H(y \vee H^{-1}(\widehat{W}^{(n)}(t))) - H(y \vee \widehat{F}^{(n)}(t))]$$
$$\Rightarrow J^*(y \vee F^*(t)) - J^*(F^*(t))\mathbb{I}_{\{y < F^*(t)\}}$$
$$= J^*(y)\mathbb{I}_{\{F^*(t) \leq y\}}$$

in the $M_1$ topology. □

PROOF OF PROPOSITION 3.5. Recall from (6.4) that

(6.6)
$$n^{1/4}[\widehat{W}^{(n)}(t) - H(\widehat{F}^{(n)}(t))] \Rightarrow J^*(F^*(t)).$$

We derive (3.7) from (6.6) using the delta method. Let us observe that $H$ is a convex, decreasing function on $\mathbb{R}$. In particular, $H^{-1}$ is convex on $[0, \infty)$ and both $H$ and $H^{-1}$ have one-sided derivatives at each point of their domains. Clearly, $H'(y+) = G(y) - 1$, $y \in \mathbb{R}$. Also,

$$(6.7) \qquad (H^{-1})'(H(y)-)H'(y+) = 1, \qquad y < y^*.$$

Indeed, for every $\varepsilon > 0$ and $y < y^*$, we have

$$1 = \frac{1}{\varepsilon}(H^{-1}(H(y + \varepsilon)) - H^{-1}(H(y)))$$

$$= \frac{H^{-1}(H(y + \varepsilon)) - H^{-1}(H(y))}{H(y + \varepsilon) - H(y)} \cdot \frac{H(y + \varepsilon) - H(y)}{\varepsilon}$$

$$\to (H^{-1})'(H(y)-)H'(y+)$$

as $\varepsilon \downarrow 0$, so (6.7) holds. For every $n$, there is an element $D^{(n)}$ (depending on the elementary event $\omega$) belonging to the subdifferential of $H^{-1}$ at some intermediate point $z \in [\widehat{W}^{(n)}(t) \wedge H(\widehat{F}^{(n)}(t)), \widehat{W}^{(n)}(t) \vee H(\widehat{F}^{(n)}(t))]$ such that

$$(6.8) \qquad H^{-1}(\widehat{W}^{(n)}(t)) - \widehat{F}^{(n)}(t) = D^{(n)}[\widehat{W}^{(n)}(t) - H(\widehat{F}^{(n)}(t))].$$

By (6.7), for any $z > 0$,

$$(6.9) \qquad (H^{-1})'(z-) = \frac{1}{H'(H^{-1}(z)+)} = \frac{1}{G(H^{-1}(z)) - 1}.$$

Applying (6.9) to a sequence $z_n \downarrow z$, we get

$$(6.10) \qquad (H^{-1})'(z+) = \frac{1}{G(H^{-1}(z)-) - 1}.$$

Using (6.9), (6.10) and convexity of $H^{-1}$, we get

$$(6.11) \qquad \begin{aligned} &\frac{1}{G(H^{-1}(\widehat{W}^{(n)}(t) \wedge H(\widehat{F}^{(n)}(t)))) - 1} \\ &\qquad \le D^{(n)} \le \frac{1}{G(H^{-1}(\widehat{W}^{(n)}(t) \vee H(\widehat{F}^{(n)}(t)))-) - 1}. \end{aligned}$$

Multiplying both sides of (6.8) by $n^{1/4}$, using (6.4), (6.11), together with continuity of the distribution of $F^*(t)$ and the fact that $G$ has at most countably many jumps, we obtain

$$(6.12) \qquad n^{1/4}[H^{-1}(\widehat{W}^{(n)}(t)) - \widehat{F}^{(n)}(t)] \Rightarrow \frac{J^*(F^*(t))}{G(F^*(t)) - 1}.$$

Finally, let us observe that $J^*$ is a mean-zero Gaussian process, independent of $F^*(t)$, so $(J^*, F^*(t))$ and $(-J^*, F^*(t))$ have the same distribution. Therefore, (6.12) implies (3.7).  $\square$

**7. First-in-first-out simulations.** To illustrate our results, we consider the special case of constant initial lead times, that is,

$$(7.1) \qquad G(y) = \mathbb{I}_{\{y \geq y^*\}} \qquad \text{for some } y^* \in \mathbb{R}.$$

In this case, all customers in the $n$th system arrive with initial lead time $\sqrt{n}y^*$ and EDF reduces to the well-known first-in-first-out (FIFO) service discipline. By (7.1), we have $H(y) = (y^* - y)^+$ and $H^{-1}(w) = y^* - w$ for all $w \geq 0$. Thus, by (1.6), $F^*(t) = y^* - W^*(t)$. Evaluating the right-hand side of (3.5), we get $\mathbb{E}[J^*(y_1)J^*(y_2)] = \lambda(\alpha^2\rho^2 + \beta^2)((y^* - y_1) \wedge (y^* - y_2))$, so $J^*(y) = B(y^* - y)$, $y \leq y^*$, where $B$ is a Brownian motion with zero drift and variance $\lambda(\alpha^2\rho^2 + \beta^2)$ per unit time, independent of $W^*(t)$. By Theorem 3.3,

$$(7.2) \quad n^{1/4}[\widehat{\mathcal{W}}^{(n)}(t)(y,\infty) - (y^* - (y \vee \widehat{F}^{(n)}(t)))] \Rightarrow B((y^* - y) \wedge W^*(t))$$

in $D(-\infty, y^*]$ [by definition, $\widehat{F}^{(n)}(t) \leq y^*$]. Theorem 3.4 yields

$$(7.3) \quad n^{1/4}[\widehat{\mathcal{W}}^{(n)}(t)(y,\infty) - ((y^* - y) \wedge \widehat{W}^{(n)}(t))] \Rightarrow B(y^* - y)\mathbb{I}_{\{W^*(t) \geq y^* - y\}}$$

in $D(-\infty, y^*]$ endowed with the $M_1$ topology. Finally, because $W^*(t) > 0$ almost surely, we have $G(F^*(t)) = 0$ almost surely and by Proposition 3.5,

$$(7.4) \qquad n^{1/4}[(y^* - \widehat{W}^{(n)}(t)) - \widehat{F}^{(n)}(t)] \Rightarrow B(W^*(t))$$

in $D(-\infty, y^*]$.

This case (with $y^* = 0$ for convenience) was already considered in [9]. In particular, (7.2) and (7.3) are contained in Corollary 4.5 and Theorem 4.1 from [9], respectively. Note that the convergence in Theorem 4.1 from [9] takes place in $D(-\infty, 0]$ endowed with the $M_1$ topology, although this was not written explicitly in the statement of that theorem.

We conduct two types of simulations. For the first simulation, we observe that an important use of the theory in [4] is to predict the amount of lateness that will be incurred by the queueing system, and the result of this paper concerns the accuracy of that prediction. In particular, [4] provides the limiting lateness result (see Theorem 3.2)

$$(7.5) \qquad \widehat{\mathcal{W}}^*(t)(-\infty, 0] = W^*(t) - \widehat{\mathcal{W}}^*(t)(0, \infty) = (W^*(t) - y^*)^+.$$

In particular,

$$(7.6) \qquad (\widehat{W}^{(n)}(t) - y^*)^+ - \widehat{\mathcal{W}}^{(n)}(t)(-\infty, 0] \Rightarrow 0.$$

The first term on the left-hand side of (7.6) is the lateness predicted by the theory as a function of the scaled workload, and the second term is the actual lateness. Setting $y = 0$ in (7.3), we obtain

$$(7.7) \quad n^{1/4}[(\widehat{W}^{(n)}(t) - y^*)^+ - \widehat{\mathcal{W}}^{(n)}(t)(-\infty, 0]] \Rightarrow B(y^*)\mathbb{I}_{\{W^*(t) \geq y^*\}}.$$

With the change of variable $s = nt$, we conclude from (7.7) that

$$
\begin{aligned}
(W^{(n)}(s) &- \sqrt{n}\, y^*)^+ - \mathcal{W}^{(n)}(s)(-\infty, 0] \\
&= \sqrt{n}(\widehat{W}^{(n)}(t) - y^*)^+ - \sqrt{n}\widehat{\mathcal{W}}^{(n)}(t)(-\infty, 0] \\
&\stackrel{d}{\approx} n^{1/4} B(y^*) \mathbb{I}_{\{W^*(t) \geq y^*\}} \\
&\stackrel{d}{=} B(\sqrt{n}\, y^*) \mathbb{I}_{\{\sqrt{n} W^*(t) \geq \sqrt{n} y^*\}} \\
&= B(\sqrt{n}\, y^*) \mathbb{I}_{\{\widetilde{W}(s) \geq \sqrt{n} y^*\}},
\end{aligned}
\tag{7.8}
$$

where $\widetilde{W}(s) \triangleq \sqrt{n}\, W^*(\frac{s}{n})$ is a reflected Brownian motion with drift $-\frac{\gamma}{\sqrt{n}}$ and variance $\lambda(\alpha^2 \rho^2 + \beta^2)$ per unit time.

We assume in this section that $\gamma$ is strictly positive. Under this assumption, the stationary distribution of $\widetilde{W}$ is exponential with mean $\frac{1}{\theta}$, where $\theta = \frac{2\gamma}{\sqrt{n}\lambda(\alpha^2\rho^2+\beta^2)}$. When conducting simulations on a single system, as opposed to a sequence of systems indexed by $n$, we know only the prelimit quantities $\frac{\gamma^{(n)}}{\sqrt{n}} = 1 - \rho^{(n)}$, $\lambda^{(n)}$, $\alpha^{(n)}$ and $\beta^{(n)}$ for a single value of $n$, and we do not know the value of $n$. We thus define

$$
\theta^{(n)} \triangleq \frac{2(1 - \rho^{(n)})}{\lambda^{(n)}[(\alpha^{(n)})^2(\rho^{(n)})^2 + (\beta^{(n)})^2]},
$$

a quantity we use as a surrogate for $\theta$. In particular, an approximate density for $\widetilde{W}(s)$ is

$$
f(y) = \theta^{(n)} \exp(-\theta^{(n)} y), \qquad y \geq 0,
\tag{7.9}
$$

and an approximate distribution for $(W^{(n)}(s) - \sqrt{n}\, y^*)^+ - \mathcal{W}^{(n)}(s)(-\infty, 0]$ is a mixture of a normal distribution with mean zero and variance $\frac{2(1-\rho^{(n)})\sqrt{n}\, y^*}{\theta^{(n)}}$ with total mass $\exp(-\theta^{(n)}\sqrt{n}\, y^*)$ and a point mass of size $1 - \exp(-\theta^{(n)}\sqrt{n}\, y^*)$ at 0. We will shortly present simulations to assess the accuracy of this approximate distribution.

While (7.2) and (7.3) are assertions about convergence of processes, (7.4) is an assertion about convergence of random variables. It describes the accuracy of the prediction of the frontier as the function $H^{-1}$ applied to the workload. From (7.4) we have

$$
\begin{aligned}
\sqrt{n}\, y^* - W^{(n)}(s) - F^{(n)}(s) &= \sqrt{n}\, y^* - \sqrt{n}\widehat{W}^{(n)}(t) - \sqrt{n}\widehat{F}^{(n)}(t) \\
&\stackrel{d}{\approx} n^{1/4} B(W^*(t)) \\
&\stackrel{d}{=} B(\sqrt{n}\, W^*(t)) = B(\widetilde{W}(s)).
\end{aligned}
\tag{7.10}
$$

Because $B(y)$ is normal with mean-zero and variance approximately equal to $\frac{2(1-\rho^{(n)})y}{\theta^{(n)}}$ and $B(y)$ is independent of $W^*(t)$, the density of $B(\widetilde{W}(s))$ is approximately

$$(7.11) \quad g(x) = \sqrt{\frac{\theta^{(n)}}{4\pi(1-\rho^{(n)})}} \int_{-\infty}^{\infty} \frac{1}{\sqrt{y}} \exp\left(-\frac{\theta^{(n)}x^2}{4(1-\rho^{(n)})y}\right) f(y)\, dy.$$

Using the Laplace transform formula

$$\int_0^{\infty} \frac{1}{\sqrt{y}}\, e^{-(a/4y)-py}\, dy = \sqrt{\frac{\pi}{p}}\, e^{-\sqrt{ap}}, \qquad a \geq 0, p > 0,$$

we can simplify (7.11) to obtain Laplace's density

$$(7.12) \qquad g(x) = \frac{\theta^{(n)}}{2\sqrt{1-\rho^{(n)}}} \exp\left(-\frac{|x|\theta^{(n)}}{\sqrt{1-\rho^{(n)}}}\right), \qquad x \in \mathbb{R}.$$

To test the predictive value of (7.8) and (7.10) as approximations to the empirical distributions of $(W^n(s) - \sqrt{n}\, y^*)^+ - \mathcal{W}^{(n)}(s)(-\infty, 0]$ and $\sqrt{n}\, y^* - W^{(n)}(s) - F^{(n)}(s)$, we simulated a single server queueing system with Poisson arrivals ($\lambda^{(n)} = 0.96$) and with three different service distributions: Exponential(1), Gamma(2, 0.5) and Uniform[0.5, 1.5]. Each of these distributions has mean 1 but different variances, so in all cases, $\rho^{(n)} = 0.96$. All customers had a constant initial lead time of $\sqrt{n}\, y^* = 30$. The queueing system began in an empty state and was simulated for $T = 4{,}000$ time units to ensure that the equilibrium assumption underlying the use of (7.9) as an approximate density for $\widetilde{W}(s)$ was valid. The simulation was independently repeated a total of 4,000 times to determine the empirical distribution.

Figures 1, 2 and 3 present normal Q–Q plots of the values for $(W^{(n)}(s) - \sqrt{n}\, y^*)^+ - \mathcal{W}^{(n)}(s)(-\infty, 0]$ restricted to the situation in which $W^{(n)}(s) \geq 30$. In this case, these values should be approximated by a normal distribution and the normal plot should be close to linear. These three figures provide strong confirmation of the accuracy of the continuous part of the proposed limiting distribution. Table 1 compares the theoretical probability that $W^{(n)}(s) \geq 30$ with the empirical probability derived from the simulation. The tabled probabilities are the total mass associated with the continuous part of the distribution. Again, the table shows the theory to be remarkably accurate.

The same simulations were used to assess the second limiting result. In particular, Figures 4, 5 and 6 address the accuracy of approximating $\sqrt{n}\, y^* - W^{(n)}(s) - F^{(n)}(s)$ by the Laplace distribution defined by (7.12). The figures differ only in the particular choice of service distribution. The figures on the left present the sorted values of $\sqrt{n}\, y^* - W^{(n)}(s) - F^{(n)}(s)$ from 4,000 independent simulations of an M/G/1 queue stopped at time $T = 4{,}000$. To evaluate the appropriateness of the proposed Laplace distribution approximation, the plots on the right present a Q–Q plot with respect to the Laplace distribution. The more linear the plot, the more appropriate is the Laplace distribution. The curves in each of these figures are highly linear; hence, they offer strong evidence of the appropriateness of this approximation.

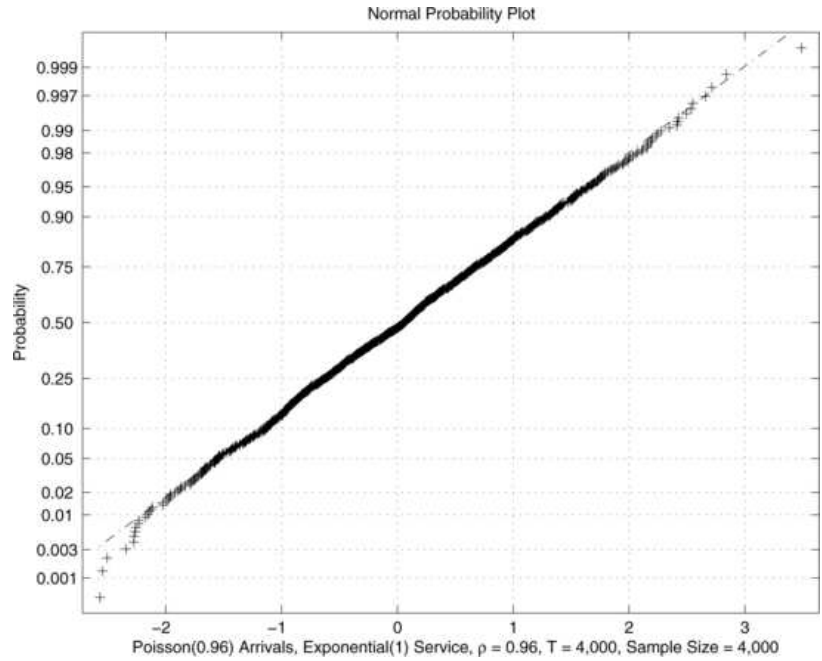

Fig. 1.   *Q–Q normal plot of continuous part of* (7.8): *exponential service.*

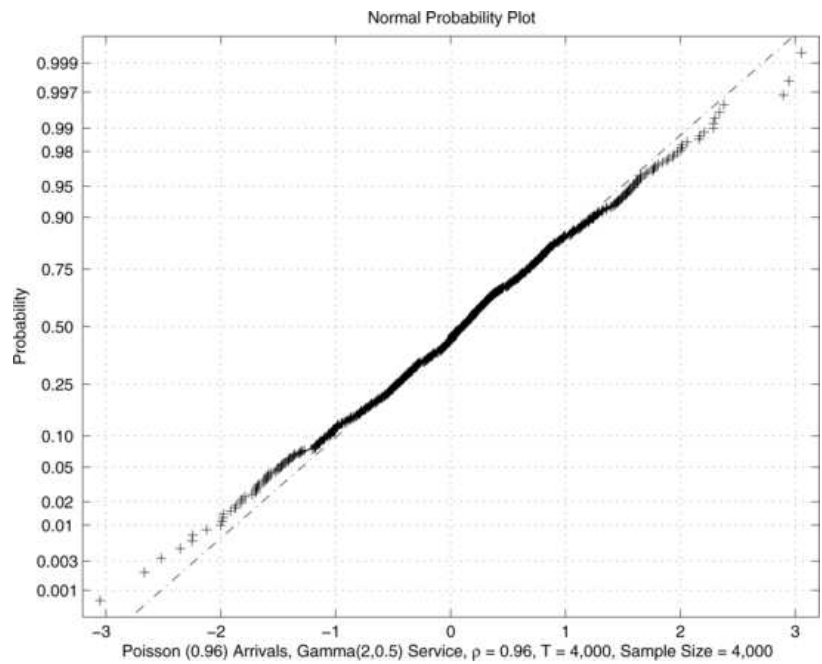

Fig. 2.   *Q–Q normal plot of continuous part of* (7.8): *gamma service.*

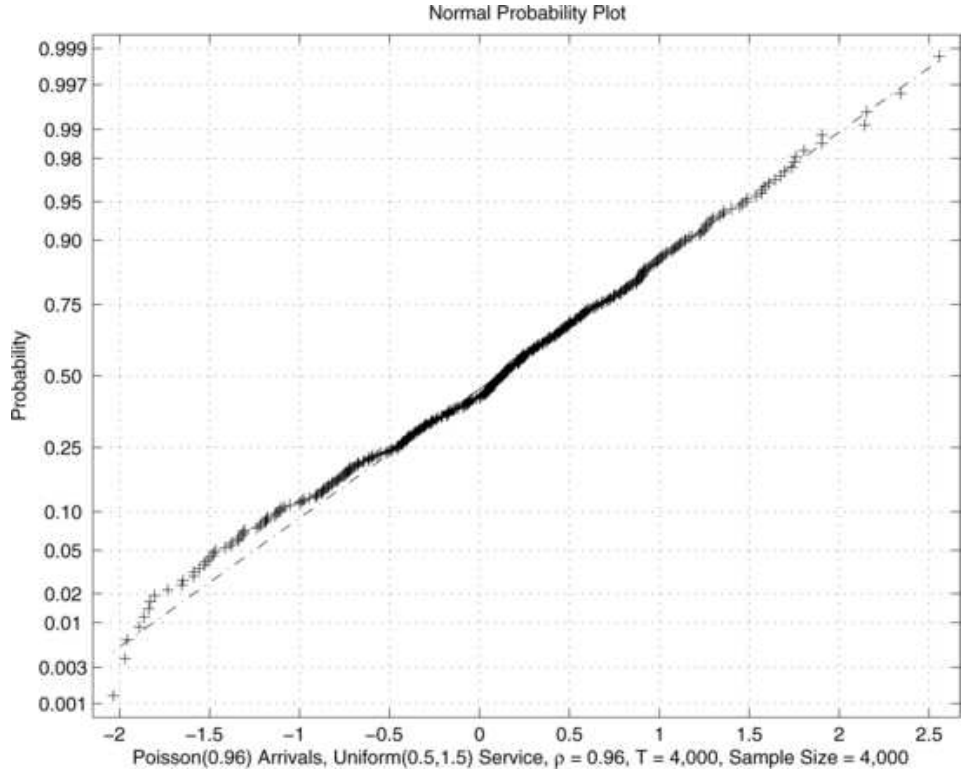

Fig. 3.   *Q–Q normal plot of continuous part of (7.8): uniform service.*

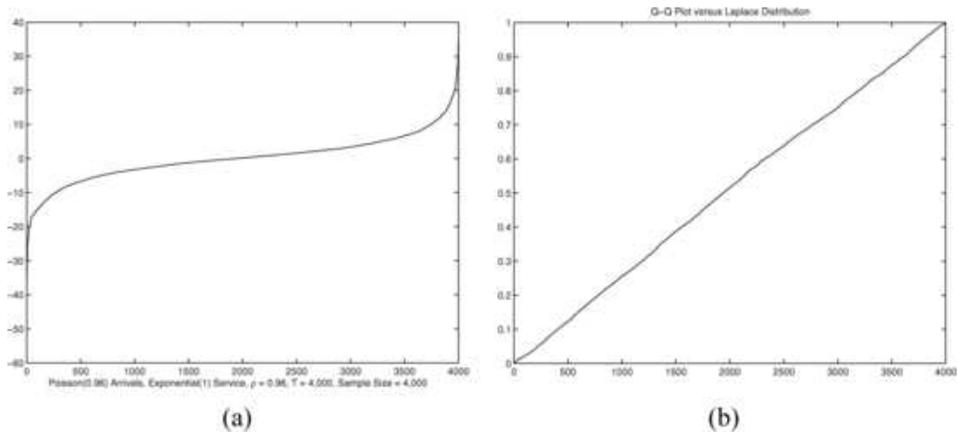

Fig. 4.   *Sorted values of (7.10) and Q–Q plot versus Laplace distribution (7.12): exponential service.*

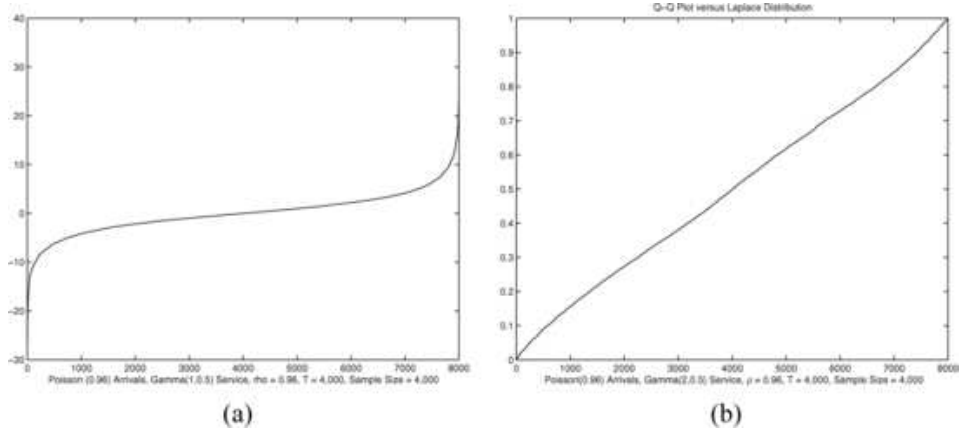

Fig. 5.  *Sorted values of (7.12) and Q–Q plot versus Laplace distribution (7.12): gamma service.*

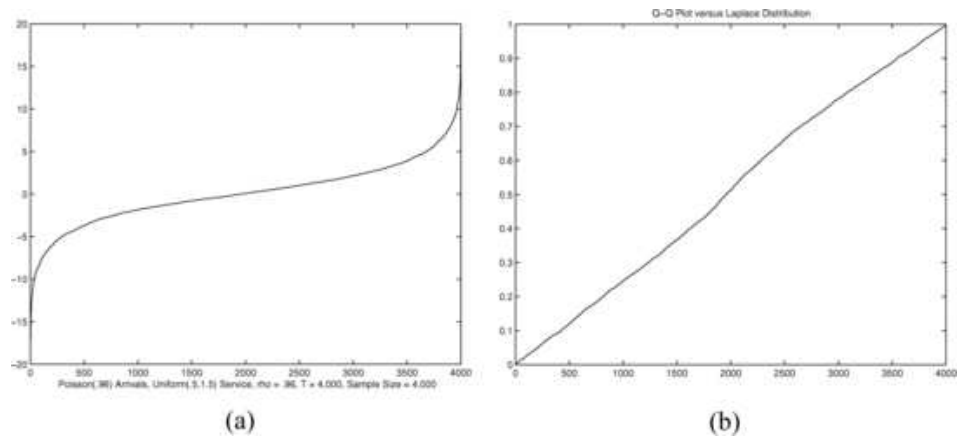

Fig. 6.  *Sorted values of (7.12) and Q–Q plot versus Laplace disribution (7.12): uniform service.*

**Acknowledgment.**  We thank Tomasz Komorowski, who pointed out that the theory of empirical processes was the right tool for analyzing the random fields appearing in this work.

TABLE 1
*Mass of the continuous part of (7.8)*

| Service distribution | Variance | Theory: $\exp(-30\theta^{(n)})$ | Simulation |
|---|---|---|---|
| Exponential | 1.000 | 0.2865 | 0.2840 |
| Gamma$(2, 0.5)$ | 0.500 | 0.1889 | 0.1880 |
| Uniform | 0.083 | 0.0995 | 0.0975 |

## REFERENCES

[1] BICKEL, P. J. and WICHURA, M. J. (1971). Convergence for multiparameter stochastic processes and some applications. *Ann. Math. Statist.* **42** 1656–1670. MR0383482

[2] BILLINGSLEY, P. (1986). *Probability and Measure*, 2nd ed. Wiley, New York. MR0830424

[3] BILLINGSLEY, P. (1999). *Convergence of Probability Measures*, 2nd ed. Wiley, New York. MR1700749

[4] DOYTCHINOV, B., LEHOCZKY, J. P. and SHREVE, S. E. (2001). Real-time queues in heavy traffic with earliest-deadline-first queue discipline. *Ann. Appl. Probab.* **11** 332–378. MR1843049

[5] ETHIER, S. N. and KURTZ, T. G. (1985). *Markov Processes*: *Characterization and Convergence.* Wiley, New York. MR0838085

[6] IGLEHART, D. and WHITT, W. (1970). Multiple channel queues in heavy traffic. I. *Adv. in Appl. Probab.* **2** 150–177. MR0266331

[7] IGLEHART, D. and WHITT, W. (1971). The equivalence of functional central limit theorems for counting processes and associated partial sums. *Ann. Math. Statist.* **42** 1372–1378. MR0310941

[8] KARATZAS, I. and SHREVE, S. E. (1988). *Brownian Motion and Stochastic Calculus.* Springer, New York. MR0917065

[9] KRUK, Ł., LEHOCZKY, J. P. and SHREVE, S. E. (2003). Second-order approximation for the customer time in queue distribution under the FIFO service discipline. *Annales UMCS Informatica AI* **1** 37–48.

[10] KRUK, L., LEHOCZKY, J. P., SHREVE, S. E. and YEUNG, S.-N. (2003). Multiple-input heavy-traffic real-time queues. *Ann. Appl. Probab.* **13** 54–99. MR1951994

[11] KRUK, L., LEHOCZKY, J. P., SHREVE, S. E. and YEUNG, S.-N. (2004). Earliest-deadline-first service in heavy traffic acyclic networks. *Ann. Appl. Probab.* **14** 1306–1352. MR2071425

[12] PANWAR, S. S. and TOWSLEY, D. (1992). Optimality of the stochastic earliest deadline policy for the G/M/c queue serving customers with deadlines. *Second ORSA Telecommunications Conference.* Boca Raton.

[13] PROKHOROV, YU. (1956). Convergence of random processes and limit theorems in probability theory. *Theory Probab. Appl.* **1** 157–214. MR0084896

[14] SHORACK, G. R. and WELLNER, J. A. (1986). *Empirical Processes with Applications to Statistics.* Wiley, New York. MR0838963

[15] STANKOVIC, J. A., SPURI, M., RAMAMRITHAM, K. and BUTTAZZO, G. C. (1998). *Deadline Scheduling for Real-Time Systems.* Springer, Berlin. MR1600916

[16] WHITT, W. (2002). *Stochastic-Process Limits*: *An Introduction to Stochastic-Process Limits and Their Application to Queues.* Springer, New York. MR1876437

[17] YEUNG, S. N. and LEHOCZKY, J. P. (2002). Real-time queueing networks in heavy traffic with EDF and FIFO queue discipline. Working paper, Dept. Statistics, Carnegie Mellon Univ.

L. KRUK
INSTITUTE OF MATHEMATICS
MARIA CURIE-SKLODOWSKA UNIVERSITY
PL. MARII CURIE-SKLODOWSKIEJ 1
20-031 LUBLIN
POLAND
E-MAIL: lkruk@hektor.umcs.lublin.pl

J. LEHOCZKY
DEPARTMENT OF STATISTICS
CARNEGIE MELLON UNIVERSITY
PITTSBURGH, PENNSYLVANIA 15213
USA
E-MAIL: jpl@stat.cmu.edu

S. SHREVE
DEPARTMENT OF MATHEMATICAL SCIENCES
CARNEGIE MELLON UNIVERSITY
PITTSBURGH, PENNSYLVANIA 15213
USA
E-MAIL: shreve@andrew.cmu.edu